\documentclass[dvips,ejs]{imsart}

\RequirePackage[OT1]{fontenc}
\RequirePackage{amsthm,amsmath,natbib}
\RequirePackage{amsmath,amssymb}
\RequirePackage[colorlinks]{hyperref}
\RequirePackage{hypernat}

\pubyear{2006}
\volume{0}
\issue{0}

\startlocaldefs
\numberwithin{equation}{section}
\theoremstyle{plain}
\newtheorem{thm}{Theorem}[section]
\endlocaldefs

\newcommand{\bea}{\begin{eqnarray*}}
\newcommand{\eea}{\end{eqnarray*}}
\newcommand{\bean}{\begin{eqnarray}}
\newcommand{\eean}{\end{eqnarray}}
\newcommand{\bdm}{\begin{displaymath}}
\newcommand{\edm}{\end{displaymath}}
\newtheorem{lemma}{Lemma}

\newtheorem{corollary}{Corollary}

\def\R{\mathbb{R}}
\def\I{1\!{\rm I}}             
\def\Xn{(X_1, \cdots, X_n)}        

\def\boldfacefake #1{%
        \hbox{%
                \mathsurround=0pt
                \hbox to 0.25pt{$#1$\hss}%
                \hbox to 0.25pt{$#1$\hss}%
                \hbox {$#1$}%
        }%
}
{\hspace{-0.15 cm}
\begin{Sbox}
\begin{minipage}{12 cm}
 \vspace{0.15 cm} }%
{\vspace{0.15 cm}
\end{minipage}
\end{Sbox}
\fbox{\TheSbox}}

\newcommand{\esl}{\end{slide}}
\newcommand{\bsl}{\begin{slide}}

\newcommand{\tr}[1]{\operatorname{tr}\left(#1\right)}

\newcommand{\id}{\mathbf{id}}

\newcommand{\bX}{\mathbf{X}}

\newcommand{\bzero}
{\text{\mathversion{bold}$0$\mathversion{normal}}}

\newcommand{\bla}{\text{\mathversion{bold}$\lambda$\mathversion{normal}}}

\newcommand{\mez}{\frac 12}

\newcommand{\te}{\theta}
\newcommand{\la}{\lambda}

\newcommand{\bei}{\begin{itemize}}
\newcommand{\eni}{\end{itemize}}
\newcommand{\beq}{\begin{equation}}
\newcommand{\enq}{\end{equation}}

\newcommand {\bdot}{\hbox{\Huge .}}
\newcommand {\dotdot}{{\hbox{\Huge .}\kern-0.1667em\hbox{\Huge .}}}
\newcommand {\onedot}{1\kern-0.1667em\bdot}
\newcommand {\twodot}{2\kern-0.1667em\bdot}
\newcommand {\idot}{i\kern-0.1667em\bdot}
\newcommand {\jdot}{j\kern-0.1667em\bdot}
\newcommand {\mdot}{m\kern-0.1667em\bdot}
\newcommand {\dotj}{\kern-0.1667em\bdot\kern-0.1667em j}

\newcommand{\slopefrac}[2]{\leavevmode\kern.1em
\raise .5ex\hbox{\the\scriptfont0 #1}\kern-.1em
/\kern-.15em\lower .25ex\hbox{\the\scriptfont0 #2}}

\renewcommand{\itemsep}{\skip0}
\renewcommand{\parsep}{\skip0}
\newcounter{example}[section]
\def\theesempio{\thesection.\arabic{example}}

\def\theesempio{\thesection.\arabic{example}}

\def\R{\mathbb{R}}
\def\I{1\!{\rm l}}
\def\Xn{\underline{X}_n}        

\def\N{\mathbb{N}}

\begin{document}

\begin{frontmatter}
\title{Bayesian nonparametric estimation of the spectral density of a long memory Gaussian time series}
\runtitle{Bayesian nonparametric for long memory time series}

\begin{aug}
\author{\fnms{Judith} \snm{Rousseau}
\ead[label=e1]{rousseau@ceremade.dauphine.fr}}
\and
\author{\fnms{Brunero} \snm{Liseo}\corref{}
\ead[label=e2]{brunero.liseo@uniroma1.it}}

\address{Universit\'e Paris Dauphine, France  \\
Sapienza Universit\'a di Roma, Italy \\
\printead{e1,e2}}

\runauthor{Rousseau and Liseo}

\affiliation{Some University and Another University}

\end{aug}

\begin{abstract}
Let $\bX=\{X_t,\, t =1,2,\dots \}$ be a stationary Gaussian random process,
with mean $EX_t = \mu$ and covariance function
$\gamma(\tau) = E(X_t-\mu)(X_{t+\tau}-\mu)$.
Let  $f(\la)$ be the corresponding spectral density;
a stationary Gaussian process is said to be
long-range dependent, if the spectral density $f(\la)$ can be
written as the product of a slowly varying function $\tilde{f}(\la)$
and the quantity $\la^{-2d}$. In this paper we propose a novel
Bayesian nonparametric approach to the estimation of the spectral
density of $\bX$. We prove that, under some specific assumptions on
the prior distribution, our approach assures posterior consistency
both when $f(\cdot)$ and $d$ are the objects of interest. The rate
of convergence of the posterior sequence depends in a significant
way on the structure of the prior; we provide some general results
and also consider the fractionally exponential (FEXP) family of priors
(see below). Since it
has not a well founded justification in the long memory set-up, we
avoid using the Whittle approximation to the likelihood function and
prefer to use the true Gaussian likelihood.
\end{abstract}


\tableofcontents
\end{frontmatter}

\section{Introduction}

Let $\bX=\{X_t,\, t =1,2,\dots \}$ be a stationary Gaussian
random process, with mean $EX_t = \mu$ and covariance function
$\gamma(\tau) = E(X_t-\mu)(X_{t+\tau}-\mu)$. Let  $f(\la)$ be the corresponding spectral density, which satisfies the relation
$$
\gamma(\tau) = \int_{-\pi}^\pi f(\la) e^{it\la} d\la \qquad (\tau=0, \pm 1, \pm2,\dots) . $$
A stationary Gaussian process is said to be long-range dependent, if there exist a positive number $C$ and a value $d$ ($0<d<1/2)$ such that
$$\lim_{\la\rightarrow 0} \frac{f(\la)}{C \la^{-2d}}=1 .$$
Alternatively, one can define a long memory process as one such that its spectral density $f(\la)$ can be written as the product of a slowly varying function $\tilde{f}(\la)$ and the quantity $\la^{-2d}$ which causes the presence of a pole of $f(\la)$ at the origin.

Interest in long-range dependent time series has increased enormously
over the last fifteen years; \cite{ber:94} provides a comprehensive
introduction and the book edited by \cite{dou:03} explores in depth both theoretical aspects and various applications of long-range dependence analysis in several different disciplines, from telecommunications engineering to economics and finance, from astrophysics and geophysics to medical time series and hydrology.

Pioneering work on long memory process is due to
\cite{mb:68}, Mandelbrot and Wallis (1969) and others. Fully
parametric maximum likelihood estimates of $d$ were
introduced in the Gaussian case by \cite{ft:86} and
\cite{dhl:89} and they have recently been developed in much
greater generality by \cite{gt:99}; a regression approach
to the estimation of the spectral density of long memory time
series is provided in \cite{gewe:83}; generalised linear
regression estimates were suggested by \cite{ber:93}. However, parametric inference can be highly biased under mis-specification of the true model: this fact has suggested semiparametric approaches: see for
instance \cite{rob:95a}.

\noindent
Due to factorization of the spectral density $f(\la)=\la^{-2d}\,\tilde{f}(\la)$, a semiparametric approach to inference seems particularly appealing in this context. One needs to estimate $d$ as a measure of long-range dependence while no particular modeling assumptions on the structure of the covariance function at short ranges are necessary: \cite{lmp:01} consider a Bayesian approach for this problem, while \cite{bardet:03} provides an exhaustive review on the classical approaches.

\noindent
Practically all the existing procedures either exploit the regression structure of the log-spectral density in a reasonably small neighborhood of the origin \cite{rob:95a} or use an approximate likelihood function based on the so called Whittle's approximation \cite{whi:62}, where the original data vector $\Xn=(X_1, X_2, \dots, X_n)$ gets transformed into the periodogram $I(\la)$ computed at the Fourier frequencies $\la_j= 2\pi\,j/n, \,\,j=1,2,\dots,n$,
and the ``new" observations $I(\la_1),\dots,I(\la_n)$ are, under a short range dependence, approximately independent, each $I(\la_j)/f(\la_j)$ having an exponential distribution. This is for example the approach taken in \cite{chgh:04}, which develop a Bayesian nonparametric analysis for the spectral density of a short memory time series. Unfortunately, the Whittle's approximation fails to hold in the presence of long range dependence, at least for the smallest Fourier frequencies.

In this paper we propose a Bayesian nonparametric approach to the
estimation of the spectral density of the stationary Gaussian
process: we avoid the use of the Whittle approximation and we deal
with the true Gaussian likelihood function.

The literature on Bayesian nonparametric inference has increased
tremendously in the last decades, both from a theoretical and a
practical point of view. Much of this literature has dealt with
the independent case, mostly when the observations are identically
distributed. The theoretical perspective was mainly dedicated to
either construction of processes used to define the prior
distribution with finite distance properties of the posterior, in
particular when such a prior is conjugate, see for instance
\cite{ghrm:00} for a review on this, or to consistency and
rates of convergence properties of the posterior, see for
instance \cite{ggvdv:01} or \cite{shw:01}.

The dependent case has hardly been considered from a theoretical
perspective apart from \cite{chgh:04}, who deal with Gaussian
weakly dependent data and, in a more general setting, \cite{gvdv:06}.
In this paper we study the asymptotic
properties of the posterior distributions for Gaussian
long-memory processes, where the unknown parameters are the
spectral density and the long-memory parameter $d$. General
consistency results are given and a special type of prior, namely
the FEXP prior as it is based on the FEXP model, is studied. From
this, consistency of  Bayesian estimators of both the spectral
density and the long memory parameter are obtained.  To
understand better the link between the Bayesian and the
frequentist approaches we also study the rates of convergence of
the posterior distributions, first in a general setup and then in
the special case of FEXP priors. The approach considered here is
similar to what is often used in the independent and identically
distributed case, see for instance \cite{ggvdv:01}. In
particular we need to control prior probability on some
neighborhood of the true spectral density and to control a sort
of entropy of the prior (see Section \ref{sec:rates}); however the
techniques are quite different due to the dependence structure of
the process.

The gist of the paper is to provide a fully nonparametric Bayesian analysis of long range dependence models.
In this context there already exist many elegant and maybe more general (in the sense of being valid even without the
Gaussian assumption) classical solutions. However we believe that a Bayesian solution would be still important because of
the following reasons.
\bei
\item[i)] By definition, our scheme allows to include in the analysis some
prior information which may be available in some applications.
\item[ii)] While classical solutions are, in a way or another, based on some asymptotic arguments, our Bayesian approach
relies only on the (finite sample size) observed likelihood function (and prior information).
\item[iii)] We are able to provide a valid approximation to the ``true" posterior distribution of the main parameters of
 interest in the model, namely the long memory parameter $d$ or the global spectral density.
\eni Also, on a more theoretical perspective, we believe that this
paper can be useful to clarify the intertwines between Bayesian
 and frequentist approaches to the problem.

The paper is organized as follows: in the next section we first
introduce the necessary notation and mathematical objects; then we
provide a general theorem which states some sufficient condition
to ensure consistency of the posterior distribution. We also
discuss in detail a specific class of priors, the FEXP prior,
which takes its name after the fractional exponential model which
has been introduced by
\cite{robi:91} (see also \cite{robi:94} to model the spectral
density of a covariance stationary long-range dependent process. The FEXP model can be seen as a generalization of the exponential model proposed by \cite{blom:73} and it allows for semi-parametric modeling of long range dependence;
see also \cite{ber:94} or \cite{hms:02}.
In Section
\ref{sec:rates} we study the rate of convergence of the posterior
distribution  first in the general case and then in the case of
FEXP priors. The final section is devoted to discussion of related problems.

\section{Consistency results} \label{sec:consG}

We observe a set of $n$ consecutive realizations $\Xn=(X_1, \dots,
X_n)$ from a Gaussian stationary process with spectral density  $f_0$, where $f_0(\lambda) =
|\lambda|^{-2d_0}\tilde{f}_0(\lambda)$.
Because of the Gaussian assumption, the density of $\Xn$ can be written as
\begin{eqnarray}
\label{dens:gaus}
\varphi_{f_0}(\Xn) = \frac{ e^{-\Xn'T_n(f_0)^{-1}\Xn /2 } }
{|T_n(f_0)|^{1/2}(2\pi)^{n/2}}, \end{eqnarray}
where $T_n(f_0)= [\gamma(j-k)]_{1\leq j,k\leq n}$ is the covariance matrix with a Toeplitz structure. The aim is to
estimate both $\tilde{f}_0$ and $d_0$ using Bayesian nonparametric
methods.

Let ${\cal F} = \{ f, f \,\mbox{symmetric on }
[-\pi,\pi], \int |f| < \infty\}$ and ${\cal F}_+ = \{ f \in {\cal
F}, f \geq 0 \}$; then ${\cal F}_+$ denotes the set of spectral
densities. We first define three types of pseudo-distances on ${\cal
F}_+$. The Kullback-Leibler divergence for finite $n$ is defined as
 \begin{eqnarray*}
 KL_n(f_0;f) &=& \frac{ 1 }{ n }\int_{\R^n }\varphi_{f_0}(\Xn) \left[
  \log{\varphi_{f_0}}(\Xn) -   \log{\varphi_{f}}(\Xn) \right] d\Xn \\
   &=&
 \frac{ 1 }{ 2n } \left\{ \tr{T_n(f_0)T_n^{-1}(f) - \id } - \log \det(T_n(f_0)T_n^{-1}(f))\right\}
\end{eqnarray*}
where $\id$ represents the identity matrix of the appropriate order.
Letting  $n\rightarrow \infty$, we can define, when it exists, the
quantity
 $$
 KL_\infty(f_0;f) = \frac{ 1 }{ \pi }\int_{-\pi}^\pi \left [
 \frac{ f_0(\la) }{ f(\la)} - 1 - \log{ \frac{ f_0(\lambda) }{ f(\lambda) } } \right ]d\lambda .$$
We also define two symmetrized version of $KL_n$, namely
$$h_n(f_0, f)=
KL_n(f_0;f)+KL_n(f;f_0); \,
d_n(f_0,f)
= \min \{ KL_n(f_0;f),KL_n(f;f_0)\}$$
and their corresponding limits as $n\rightarrow \infty$:
$$
h(f_0,f) = \frac{ 1 }{ 2\pi }\int_{-\pi}^\pi\left [ \frac{f_0(\lambda)}{f(\lambda)}
 + \frac{ f(\lambda) }{ f_0(\lambda) } -2 \right ] d\lambda;$$
$$d(f_0,f) = \min \{ KL_\infty (f_0;f), KL_\infty(f;f_0) \}.$$
We also consider the $\mbox{L}_2$ distance between the logarithms
of the spectral densities, namely
 \begin{eqnarray}\label{def:L2}
 \ell(f,f') &=& \int_{-\pi}^\pi (\log{f}(\lambda) -
 \log{f'}(\lambda))^2d\lambda.
 \end{eqnarray}
This distance has been considered in particular by
\cite{mousou:03}. This is quite a natural distance in the
sense that it always exists, whereas the $\mbox{L}_2$ distance
between $f$ and $f'$ need not, at least in the types of models
considered in this paper. Let $\pi$ be a prior probability
distribution on the set
$$\tilde{{\cal F}} = \{ f \in {\cal F}, f(\lambda) =
|\lambda|^{-2d}\tilde{f}(\lambda), \quad \tilde{f}\in C^0, -\mez < d
< \mez \}, \quad \tilde{{\cal F}}_+ = \{ f \in \tilde{{\cal F}}, f \geq 0
\},
$$
where $C^0$ is the set of continuous functions on $[-\pi,\pi]$.\\ Let
$A_\varepsilon = \{ f \in \tilde{ {\cal F}_+ }; d(f,f_0) \leq
\varepsilon \}$.
Our first goal is to prove the consistency of the
posterior distribution of $f_0$, that is, we show that, as $n\rightarrow \infty$,
 $$P^\pi[A_\varepsilon^c |\Xn] \rightarrow 0, \quad f_0 \quad \mbox{a.s.},$$
 where $P^\pi[.|\Xn]$ denotes the posterior distribution associated
 with the prior $\pi$.
From this, we will be able to deduce the consistency of some Bayes estimators of
the spectral density $f$ and of the long memory parameter $d$.
We first state and prove the strong consistency of the posterior distribution
under very general conditions both on the prior and on the {\it
true} spectral density.
Then, building on these results, we will obtain the
consistency of a class of Bayes estimates of the spectral density,
together with the consistency of the Bayes estimates of the long
memory parameter $d$.
The already introduced FEXP class of prior will be then proposed, and its use will be explored in detail.

\subsection{The main result} \label{subsec:main}
In this section we derive the main result about consistency of the
posterior distribution. We also discuss the asymptotic behavior of
the posterior point estimates of some parameter of major interest,
such as the long memory parameter $d$ and the global spectral
density.

\noindent Consider the following two subsets of ${\cal F}$
\begin{eqnarray}
\label{set:unif0}
 &\,& {\cal G}(d,M,m,L,\rho) = \\ \notag
 &\,&
 \left\{ f \in \tilde{{\cal F}}_+ ; f(\lambda)= |\lambda|^{-2d}\tilde{f}(\lambda), m
 \leq \tilde{f}(\lambda) \leq M, \left|\tilde{f}(x) - \tilde{f}(y)\right| \leq L|x-y|^\rho \right\},
\end{eqnarray}
where
 $-1/2  < d < 1/2$, $m, M, \rho >0$;
 \begin{eqnarray}
 \label{set:unif1}
 &\,& {\cal F}(d,M,L,\rho) = \\ \notag
 &\,& \{ f \in \tilde{{\cal F}} ; f(\lambda) = |\lambda|^{-2d}\tilde{f}(\lambda), |\tilde{f}(\lambda)| \leq
 M, \left|\tilde{f}(x) - \tilde{f}(y)\right| \leq L|x-y|^\rho\}.
\end{eqnarray}

\noindent The boundedness constraint on $\tilde{f}$ in the
definition of $\mathcal{G}(d,M,m,L,\rho)$ is here to guarantee the
identifiability of $d$, while the Lipschitz-type condition on
$\tilde{f}$, in both definitions, are actually needed to ensure
that normalized traces of products of Toeplitz matrices, that
typically appear in the distances considered previously, will
converge. We also consider the following set of spectral
densities, which is of interest in the study of rates of
convergence: let
$$ \mathcal{L}^\star(M,m,L) = \{ h(\cdot) \geq 0, \,0< m \leq h(\cdot) \leq M,\, |h(x)
- h(y)| \leq L|x-y| (|x|\wedge |y|)^{-1} \}$$ and
$$ \mathcal {L} (d,M,m,L) = \{ f=|\lambda|^{-2d} \tilde{f}(\lambda), \tilde{f} \in \mathcal{L}^\star(M,m,L) \}.$$
Note that $\mathcal{G}$ and $\mathcal{L}$ are similar, with only a
slight modification on the Lipschitz condition. The set
$\mathcal{L}$ has been considered in particular in
\cite{mousou:03}.

\noindent We now consider the main result on the consistency of the
posterior distribution. Let
$$\bar{{\cal G}}(t,M,m,L,\rho) =
\cup_{-1/2+t \leq d \leq 1/2-t}{\cal G}(d,M,m,L,\rho)$$
and
$$\bar{
\mathcal {L}} (t,M,m,L) =  \cup_{-1/2+t \leq d \leq 1/2-t}\mathcal
{L} (d,M,m,L).$$
In the following theorem and in its proof we
consider spectral densities in sets either in the form
$\bar{\mathcal G}$ or in the form $\bar{\mathcal L}$. To simplify
the presentation we give results for densities in $\bar{\mathcal G}$
only, however the results remain valid for densities in
$\bar{\mathcal L}$, the only difference being that in the conditions
in the form $4|d-d_0| \leq \gamma$ where $\gamma = \rho\wedge \rho_0
\wedge 1/2$ the quantities $\rho, \rho_0$ can be chosen equal to 1
if the corresponding spectral densities belong to $\bar{\mathcal
L}$.

\begin{thm}\label{th:cons}
Assume that there exist $(t_0,M_0,m_0,L_0)$ such that we have
\bei
\item either $f_0 \in \bar{{\cal G}}(t_0,M_0,m_0,L_0,\rho_0) $ with $0<
 \rho_0 \leq 1$
 \item or $f_0 \in \bar{ \mathcal {L}} (t_0,M_0,m_0,L_0)=
\cup_{-1/2+t_0 \leq d \leq 1/2-t_0}\mathcal {L} (d,M_0,m_0,L_0)$.
\eni
Let $t,m,M,L$ be positive reals with $t <  (\rho\wedge \rho_0)/4$.
Let $\pi$ be a prior distribution on either $\bar{{\cal
G}}(t,M,m,L,\rho)$, $\rho>0$ or $\bar{
 \mathcal L}(t,m,M,L)$.
  If the prior satisfies:
\begin{itemize}
\item[i)] $\forall \varepsilon >0$, $\pi({\cal B}_\varepsilon
)>0$, where
 $$ {\cal B}_\varepsilon = \left \{ f \in \bar{{\cal F}}_+(t,M,m) :
 h(f_0,f) \leq \epsilon, 8(d_0-d) < \gamma - t \right \}.$$

\item[ii)] For all $\varepsilon > 0$, small enough, there exists ${\cal F}_n
\subset\tilde{{\cal F}}_+ $, such that $\pi({\cal F}_n^c) \leq
e^{-nr}$
 and a smallest possible net ${\cal H}_n \subset \bar{\mathcal G}(t,M,m,L,\rho)
 ; d(f,d_0)> \varepsilon/2\}$ (resp.
 $\bar{\mathcal L}(t,M,m,L)$) such that when $n$ is
large enough, $\forall f \in {\cal F}_n \cap A_\varepsilon^c, \quad
\exists f_i \in {\cal H}_n$, $ f_i  = |x|^{-2d_i} \tilde{f}_i(x)
\geq f $ such that $4(d_i-d) \leq \rho\wedge 1/2$.

{\bf -} If $4|d-d_0| \leq \gamma -t$,  $$\frac{ 1 }{ 2 \pi }
\int_{-\pi}^{\pi}\frac{ (f_i-f)(x) }{ f_0(x) }dx \leq
h(f_0,f_i)/4$$
{\bf -} If $4(d-d_0) > \gamma -t$
$$\frac{ 1 }{ 2 \pi }
\int_{-\pi}^{\pi}\frac{ (f_i-f)^2(x) }{ f^2(x) }dx \leq b(f_0,f_i)
 |\log{\varepsilon}|^{-1}$$
{\bf -} If $4(d_0-d) > \gamma -t$
$$\frac{ 1 }{ 2 \pi }
\int_{-\pi}^{\pi}\frac{ (f_i-f)^2(x) }{ f_0^2(x) }dx \leq
b(f_i,f_0), \quad \frac{ 1 }{ 2 \pi } \int_{-\pi}^{\pi}\frac{
(f_i-f)(x) }{ f(x) }dx \leq b(f_i,f_0) |\log{\varepsilon}|^{-1}$$

\noindent
Denote by $N_n$ the logarithm of the cardinality of the smallest
possible net ${\cal H}_n$. Then, if
 $$ N_n \leq nc_1, \quad \mbox{with} \quad  c_1 \geq  \varepsilon  |\log{\varepsilon}|^{-2} ,
 $$
  then
\begin{eqnarray}\label{asconv}
P^\pi\left[ A_\varepsilon |\Xn \right] \rightarrow 1, \quad  f_0
\quad \mbox{a.s.}
\end{eqnarray}
 \end{itemize}
\end{thm}

\begin{proof}
See Appendix B.
\end{proof}

\noindent The above theorem is important to clarify which conditions
on the prior distribution $\pi$ are really crucial in a long memory
setting, where the techniques usually adopted in the i..i.d.\,case,
cannot be used and even the adoption of a Whittle approximation is
not legitimate in this setting (at least at the lowest frequencies).
From a practical perspective, however, the hardest part of the program
is actually to verify whether a specific type of priors actually
meets the conditions listed in Theorem \ref{th:cons}.
 We will discuss in detail these issues in the
 context of the FEXP prior in \S\ref{sec:cons:fexp}.

\subsection{Consistency of estimates for some quantities of interest}
We now discuss the problem of consistency for the Bayes estimates of
the spectral density. The usual quadratic loss function for the class of functions ${\cal F}$ is not the
natural one for this problem, since there exist some
spectral densities in ${\cal F}$ that are not square integrable (i.e. if
$d > 1/4$). A more reasonable loss function is the quadratic
loss on the logarithm of $f$,  as defined by (\ref{def:L2}), which
is always integrable, at least in the framework considered in the
paper. The Bayes estimator of $f$ associated with the loss $\ell$
and the prior $\pi$ is given by
$$
\hat{f}(\lambda) = \exp\{ E^\pi[ \log{f}(\lambda)| \Xn ] \} = |\lambda|^{-2\hat{d}}
 \exp\{ E^\pi[ \log{\tilde{f}}(\lambda)| \Xn ] ,
$$
where $\hat{d}=E^\pi[d|\Xn]$. Note also that the Bayes estimator of
$f$ associated with the loss $h(.,.)$ and the prior $\pi$ is given
by
\beq
\label{loss}
\hat{f}_2^\pi(\lambda) = \sqrt{ \frac{ E^\pi[f(\lambda)|\Xn] }{
 E^\pi[f^{-1}(\lambda)|\Xn] }}.
\enq
   Also, in many
applications of long memory processes, the real parameter of interest is just $d$, the long
memory exponent. It is possible to deduce,
 from Theorem \ref{th:cons}, that the  posterior mean of $d$, that is the Bayes estimator associated
 with the quadratic loss on $d$, is actually consistent.
 Let $$\bar{{\cal F}}_+(t,M,m) = \cup_{-1/2+t \leq d \leq 1/2-t}\{ f \in {\cal F}_+,
f(\lambda)= |\lambda|^{-2d}\tilde{f}(\lambda), 0<m \leq \tilde{f}
\leq M \}.$$

\begin{corollary} \label{Cor:cons:d}
Under the assumptions of Theorem \ref{th:cons},  for all
$\epsilon>0$, as $n\rightarrow \infty$,
 $$\pi\left[\{ f = |\lambda|^{-2d} \tilde{f}; |d-d_0|>
 \epsilon\}|\Xn\right] \rightarrow 0 \quad  f_0 \mbox{ a.s}$$
and $ \hat{d}  \rightarrow d_0, \quad f_0 \mbox{ a.s}.$
\end{corollary}
\begin{proof}
The result comes from the fact that, when $|d-d_0|> \epsilon$, there
exists a positive constant $\epsilon'$ such that for all $f,f_0 \in
\bar{{\cal F}}_+(t,M,m),$ $f= |\lambda|^{-2d}\tilde{f}$ and $f_0 =
|\lambda|^{-2d_0}\tilde{f}_0$, $h(f,f_0)> \epsilon'$.
In fact, assume without
loss of generality, that $d>d_0$, then for all $A > 4M/m$
\begin{eqnarray*}
h(f,f_0) &=& \frac{ 1 }{ 2\pi }\int_{-\pi}^\pi (f/f_0+f_0/f -
2)(\lambda d\lambda \\
 &\geq& \frac{ m }{ 8\pi
 M}\int_{|\lambda|^{-2(d-d_0)}>A}|\lambda|^{-2(d-d_0)}d\lambda .
\end{eqnarray*}
The above quantity is infinite if $2(d-d_0) \geq 1$, otherwise
\begin{eqnarray*}
h(f,f_0) & \geq & \frac{ m }{ 4\pi
 M (1 -2(d-d_0))}A^{1-1/(2(d-d_0))} \\
  &\geq&  \frac{ m }{ 4\pi
 M (1 -2(d-d_0))}A^{1 -1/2\epsilon} = \epsilon'.
\end{eqnarray*}
 This implies that
 $$\pi[A_{\epsilon'}^c|X] \geq \pi\left[\{ f = |\lambda|^{-2d} \tilde{f}; |d-d_0|>
 \epsilon\}|\Xn\right] \rightarrow 0, \quad f_0 \mbox{ a.s.}$$
 Since $d$ is bounded, a simple application of the Jensen's inequality gives
 \begin{eqnarray*}
 (\hat{d} - d_0)^2 \leq E^\pi[(d-d_0)^2|\Xn]  \rightarrow 0, \quad
 f_0 \mbox{ a.s.}
\end{eqnarray*}
\end{proof}
It is also possible to derive consistency results for the
\textit{point} estimate of the whole spectral density:
\begin{corollary} \label{Cor:cons:f:1}
Under the assumptions of Theorem \ref{th:cons}, if $\hat{f}_2^\pi$ is as defined in (\ref{loss}), as $n\rightarrow
\infty$,
$$h(f_0,\hat{f}_2^\pi)\rightarrow 0, \quad f_0 \mbox{ a.s.}$$
\end{corollary}

\begin{proof}
To simplify the notations we set $C$ to be a generic positive
constant. Let $H(x ) = x + x^{-1}-2$, then for any $a>0$
\begin{eqnarray*}
h(f_0,\hat{f}_2^\pi) &=& C \int\left( \frac{ \sqrt{
E^\pi[f/f_0(\lambda)|\Xn]} }{ \sqrt{ E^\pi[f_0/f(\lambda)|\Xn]} }+
\frac{ \sqrt{ E^\pi[f_0/f(\lambda)|\Xn]} }{ \sqrt{
E^\pi[f/f_0(\lambda)|\Xn]}} - 2 \right)d\lambda \\
 &\leq&
 C\int_{\lambda > a}E^\pi[ H(f/f_0(\lambda))|\Xn] +
 C\int_{\lambda<a} H\left[\frac{ \sqrt{
E^\pi[f/f_0(\lambda)|\Xn]} }{ \sqrt{ E^\pi[f_0/f(\lambda)|\Xn]}
}\right]d\lambda \\
 &=& I_1 + I_2
\end{eqnarray*}
We have:
\begin{eqnarray*}
I_1 &\leq& C\epsilon + C E^\pi \left[\I_{h(f,f_0)> \epsilon}
 \int_{\lambda > a}H(f/f_0(\lambda))d\lambda |\Xn\right]
\end{eqnarray*}
Now, consider the test $\phi_n$ defined in the proof of Theorem
\ref{th:cons} and the same type of inequality as those used in the
proof of the same theorem: for all $f \in A_\epsilon^c$
 $$E_f[1-\phi_n] \leq e^{-n\epsilon |\log{\epsilon}|^{-1}}.$$
 Then we choose a small $\delta>0$ such that
\begin{eqnarray*}
P_0^n\left[ I_1> 2C\epsilon \right] &\leq& E_0^n\left[\phi_n\right]
+ \frac{ C }{ n^3 } +
 \frac{ Ce^{n\delta} }{ \epsilon } \int_{\lambda>a} \int
 H(f/f_0(\lambda)) E_f^n\left[1 - \phi_n\right]d\pi(f)d\lambda \\
  &\leq& C n^{-3} + Ce^{-n\epsilon |\log{\epsilon}|^{-1}/2}\left[ 1 + \int_{\lambda>a}
 H(f/f_0(\lambda)) d\pi(f)d\lambda \right]\\
 &\leq & C n^{-3} + Ce^{-n \epsilon |\log{\epsilon}|^{-1}/2} a^{-2}.
\end{eqnarray*}
Let $a = \exp( -n\epsilon |\log{\epsilon}|^{-1}/ 8)$ then
\begin{eqnarray*}
P_0^n\left[ I_1> 2C\epsilon \right] &\leq&  Ce^{-n\epsilon
|\log{\epsilon}|^{-1}/4}.
\end{eqnarray*}
We also have
\begin{eqnarray*}
I_2 &=& \int_0^{a} H\left[\frac{ \sqrt{ E^\pi[f/f_0(\lambda)|\Xn]}
}{ \sqrt{ E^\pi[f_0/f(\lambda)|\Xn]}
}\right]d\lambda \\
 &\leq&
 C \int_{0<\lambda<a} \frac{ \sqrt{
E^\pi[\lambda^{-2(d-d_0)}|\Xn]} }{ \sqrt{
E^\pi[\lambda^{-2(d_0-d)}|\Xn]} } +
 \frac{ \sqrt{
E^\pi[\lambda^{-2(d_0-d)}|\Xn]} }{ \sqrt{
E^\pi[\lambda^{-2(d-d_0)}|\Xn]} } d\lambda \\
&\leq& C \int_{0<\lambda<a} \lambda^{2d_0} \lambda^{-\hat{d}}
\sqrt{E^\pi[\lambda^{-2d}|\Xn]}d\lambda + C \int_{0<\lambda<a}
\lambda^{-2d_0} \lambda^{\hat{d}}
\sqrt{E^\pi[\lambda^{2d}|\Xn]}d\lambda
\end{eqnarray*}
Let $A$ be the set where $\hat{d} = E^\pi[d|\Xn]$ converges to
$d_0$; then $P_0^\infty[A]=1$ and  $\forall \delta>0$ and $n$ large
enough,
\begin{eqnarray*}
I_2&\leq& C \int_{0<\lambda<a} \lambda^{d_0-\delta-1/2+t} d\lambda
 + C \int_{0<\lambda<a}
\lambda^{-d_0-1/2-\delta}d\lambda \\
&\leq& C(a^{1/2+d_0 -\delta}+a^{1/2-d_0 -\delta}) \\
 &\leq & e^{-nc\epsilon |\log{\epsilon}|^{-1}},
\end{eqnarray*}
for some $c>0$.
\end{proof}
\begin{corollary} \label{Cor:cons:f:2}
Under the assumptions of Theorem \ref{th:cons}, as $n\rightarrow \infty$,
$$\ell(f_0,\hat{f}^\pi) \rightarrow 0, \quad f_0 \mbox{ a.s.}$$
\end{corollary}
\begin{proof}
Note that for all $x \in \R$, $e^{x}+e^{-x} - 2 \geq x^2$.
Then
$h(f,f_0) \geq l(f,f_0)$ and
 $$P^\pi[ f; l(f,f_0) > \epsilon|\Xn] \leq P^\pi[A_\epsilon^c|\Xn]$$
 This implies, together with the fact that $l(f,f_0)$ is bounded
 when $f \in \mathcal G(t,m,M,\rho)$, that $\forall \epsilon>0$,
 $$l(\hat{f},f_0) \leq E^\pi[l(f,f_0)|\Xn] \leq \epsilon +
 CP^\pi[A_\epsilon^c|\Xn].$$
\end{proof}

Since the conditions stated in Theorem \ref{th:cons} are somewhat
non standard, they need to be carefully checked for the specific class of priors
 one is dealing with. Here we consider the class of Fractionally
Exponential priors (FEXP), and we show that these priors actually
fulfill the above conditions.

\subsection{The FEXP prior} \label{sec:cons:fexp}
Consider the set of the spectral densities with the form
 $$f(\lambda) = |1 - e^{i\lambda}|^{-2d} \tilde{f}(\lambda),$$
where $\log{\tilde{f}}(\lambda) = \sum_{j=0}^K \theta_j
\cos(j\lambda)$, for some finite $k\in\N$, and assume that the true log spectral density
satisfies $\log{\tilde{f}_{0}}(\lambda) = \sum_{j=0}^\infty
\theta_{0j}\cos(j\lambda)$ (in other words, it is equal to its
Fourier series expansion), with
 $$|\tilde{f}_{0}(\lambda) - \tilde{f}_{0}(\lambda')| \leq L \frac{ |\lambda-\lambda'| }{
 |\lambda|\wedge|\lambda'|}, \quad \sum_j |\theta_{0j}|< \infty, $$
for all $\lambda$ and $\lambda'$ in $[-\pi,\pi]$. In this section our base model is presented in a slightly different way:
however it comes to the same thing since $|1 -
e^{i\lambda}|/|\lambda|$ is  continuous and strictly positive on
$[-\pi,\pi]$.

This class of densities has been considered, from a frequentist
perspective, in \cite{hms:02}. Note that there exists an
alternative and equivalent way of writing a FEXP spectral density in
which the first coefficient of the series expansion $\te_0$ is
explicitly expressed in terms of the variance of the process, that
is $\sigma^2=2\pi\,e^{\te_0}$. We will use both the
parameterizations according to notational convenience. A prior
distribution on $f$ can then be expressed as a prior on the
parameters $(d,K,\theta_0,...,\theta_K)$ in the form
$p(K)\pi(d|K)\pi(\theta|d,K)$, where $\theta =
(\theta_0,...,\theta_K)$, and $K$ represents the (random) order of
the FEXP model. Usually, $d$ is set independent of $\te$ for given
$K$  and it is also independent of $K$ itself.
 Let $\pi(d)> 0$ on $[-1/2+t, 1/2-t]$, for some $t>0$, arbitrarily small.
 Let $K$ be a priori Poisson distributed and, conditionally
 on $K$, in order to obtain a Lipschitz condition on
 $\sum_{j=1}^K\theta_j \cos{(j\lambda)}$ we consider $\theta$'s such
 that $\sum_{j=1}^k j|\theta_j| \leq B$, where $B$ is large but
 finite. This implies in particular, that the terms
 $\sum_{j=1}^K|\theta_j|$ are uniformly bounded over the supports of
 $\pi_K$.
 A possible way to formalize it, is to assume that, for given $K$, the quantity
  $S_K=\sum_j j|\theta_{j}|$ has a finite support distribution; then, setting $V_j=j|\theta_{j}|/S_K$, $j=1, \ldots,K$,
one may consider a distribution on the set $\{z \in \R^K;
z=(z_1,...,z_K), \sum z_i=1, z_i \geq 0\}$  for example:
\centerline{ $ (V_1, \ldots, V_K) \sim \mbox{Dirichlet}(\alpha_1,
\ldots, \alpha_K), $} Since the variance of the $|\theta_{j}|$'s
should be decreasing as $j$ increases, we may assume, for example, that, for all
$j$'s, $\alpha_j=O((1+j)^{-2})$. Note that if we further assume that
$S_K$ has a Gamma distribution with mean $\sum_j \alpha_j$ and
variance $\sum_j \alpha_j^2$ then we are approximately assuming
(modulo the truncation at $A$) that $\left (|\te_1|, \ldots, |\te_k|\right )$
are independent Gamma$(1,\alpha_j)$ random variables. Alternative
parameterization are also available here; for example one can assume
that $(V_1, \cdots, V_k)$ follows a logistic normal distribution
\cite{aitc:80}, which allows for a more flexible elicitation. Under
the above conditions on the prior, the posterior distribution is
strongly consistent, in terms of the distance $d(\cdot,\cdot)$, the
estimator $\hat{f}$ as described in the previous section is almost
surely consistent and so is the estimator $\hat{d}$. To prove this,
we need to show that the FEXP prior satisfies assumptions (i) and (ii).
First, we check assumption (i): let $K_\epsilon$ be such that
$\sum_{j=K_\epsilon+1}^\infty |\theta_{0j}| \leq \sqrt{\epsilon}/4,$
then $h(f_0,f_{0\epsilon}) \leq \epsilon/8$, where
 $$f_{0\epsilon} = |1 - e^{i\lambda}|^{-2d_0}\exp\left\{
 \sum_{j=0}^{K_\epsilon}\theta_{0j}\cos{j\lambda} \right\}.$$
Let $\theta = (\theta_0,...,\theta_{K_\epsilon})$ be such that
 $\sum_{j=0}^{K_\epsilon}|\theta_{0j} - \theta_j| \leq \sqrt{\epsilon}/4$
 $j=1,\ldots,K_\epsilon$ 
  If $|d-d_0| < \epsilon c$, with
  $c \leq \left( \int_{-\pi}^{\pi}|1  - e^{i\lambda}|^{-1/4}d\lambda
  \right)^{-1}/(8\pi)$,
then
 $$ h(f,f_0) \leq \frac{\epsilon}{8\pi}\int_{-\pi}^\pi
 |1-e^{i\lambda}|^{-\epsilon/4} d\lambda+
 e^{\sqrt{\epsilon}/2}\frac{\epsilon}{4} \leq \epsilon$$
 for $\epsilon>0$ small enough.
Also $\pi_{K_\epsilon}(\{ \theta :
 |\theta_j - \theta_{0j}|< \sqrt{\epsilon}/(8K_\epsilon), \forall j\leq K_\epsilon \})
 >0$, as soon as $A > \sum_j|\theta_{0j}|$. Thus assumption  (i) of Theorem \ref{th:cons} is
 satisfied.

\noindent
 Now we verify assumption (ii). Let
$\epsilon> 0$ and set
$$f_{k,d,\theta}(\la) = |1 - e^{-i\lambda}
|^{-2d} \exp\{ \sum_{j=0}^k \theta_j \cos{(j\lambda)} \},$$ where the
$\theta_j$'s satisfy the above constraint. Consider
$$\mathcal{F}_n = \{ f_{k,d,\theta}, d \in [-1/2+t,1/2-t], k \leq
k_n\},$$ where $k_n = k_0 n /\log{n}$.
Since $\pi(K \geq k_n)< e^{-nr}$,
for some $r$ depending on $k_0$, we have that
$\pi(\mathcal{F}_n^c \geq k_n)< e^{-nr}$.
Now consider spectral densities in the form,
 $$f_i(\la) = (1 - \cos{\lambda})^{-d_i} \exp\{
 -d_i \log(2)+\sum_{j=0}^k \theta_j^i \cos{j\lambda} \}.$$
 Consider
$$f(\la) = (1 - \cos{\lambda})^{-d} \exp\{
 -d \log(2)+\sum_{j=0}^k \theta_j \cos{j\lambda} \},$$
 where $d_i - c_1 \epsilon \leq d\leq d_i$, $ \theta_0^i - c_2\epsilon \leq  \theta_0 + (d_i-d)
 \log(2) \leq \theta_0^i - c_0\epsilon$, and $\sum_{j=1}^k
 |\theta_j- \theta_j^i| \leq c_0\epsilon$. Then
 \begin{eqnarray*}
 \frac{ f(\la) }{ f_i(\la) } &=& (1 - \cos{\lambda})^{d_i-d} \exp\{
 (d_i-d) \log(2)+\sum_{j=0}^k (\theta_j-\theta_j^i) \cos{j\lambda}
 \} \\
&\leq& 1
\end{eqnarray*}
and
\begin{eqnarray*}
 \frac{ f(\la) }{ f_i(\la) } &\geq &
 (1 - \cos{\lambda})^{-c_1\epsilon} e^{-(c_2+c_0)\epsilon}
 \end{eqnarray*}
Hence by choosing $c_0,c_1,c_2$ small enough, $f_i-f$ verifies the
three inequalities considered in assumption (ii) of Theorem
\ref{th:cons}. The covering number of $\mathcal F_n$ with balls
defined by the above inequalities can be bounded by
 $$ \exp(N_n  ) \leq k_n (Ck_n/\epsilon)^{k_n+2}  \leq e^{ 2k_0 n
 (-\log{\epsilon}-\log{\log{n}})}$$
 so that if $n$ is large enough
  $$N_n \leq n \epsilon|\log{\epsilon}|^{-2}$$
and assumption (ii) is satisfied.

\section{Rates of convergence} \label{sec:rates}

In this section we first provide a general theorem relating rates for
the posterior distribution to conditions on the prior. These
conditions are, in essence, similar to the conditions obtained in
the i.i.d. case; in other words there is a condition on the the
prior mass of Kullback-Leibler neighbourhoods of the true spectral
density and an entropy condition on the support of the prior. We
then present the results in the case of the FEXP prior.

\subsection{Main result}

We now present the general Theorem on convergence rates for the
posterior distribution.

\begin{thm}\label{th:rate:spde}
Let $(\rho_n)_n$ be a sequence of positive numbers decreasing to
zero, and
 ${\cal B}_n$ a ball belonging to
 $\left( \bar{{\cal G}}(t,M,m,L,\rho) \cup \bar{\mathcal{L}}(t,M,m,L) \right),$
 defined as
 $${\cal B}_n(\delta)=
   \{f(x) =
 |x|^{-2(d-d_0)}\tilde{f}(x); KL_n(f_0;f) \leq \rho_n/4, b_n(f_0,f) \leq \rho_n , |d-d_0|
 \leq \delta \},$$
 for some $\rho \in (0,1]$.
Let $\pi$ be a prior such which satisfies conditions (i) and (ii) of Theorem
\ref{th:cons}. Assume that:\\
{\rm (i).}\, There exists $\delta >0 $ such that $\pi({\cal B}_n(\delta)) \geq exp\{-n\rho_n/2\}.$ \\
{\rm (ii).}\, For all $\epsilon >0$ small enough, there exists  a
positive sequence $(\epsilon_n)_n$ decreasing to zero and
$\bar{{\cal F}}_n \subset \tilde{{\cal F}}_+ \cap \{f, d(f,f_0) \leq
\epsilon \}$, such that $\pi(\bar{{\cal F}}_n^c \cap \{f, d(f,f_0)
\leq \epsilon \}  ) \leq e^{-2 n\rho_n}$. \\
{\rm (iii).}\, Let
$$S_{n,j} = \{ f \in \bar{{\cal F}}_n;
\varepsilon_n^2j \leq h_n(f_0,f)  \leq \varepsilon_n^2 (j+1) \},$$
with $J_n \geq j \geq J_0$, with fixed $J_0>0$ and $J_n =
\lfloor\varepsilon^2 /\varepsilon_n^2\rfloor$. $\forall J_0 \leq j
\leq J_n$, there exists a smallest possible net $\bar{{\cal
H}}_{n,j} \subset S_{n,j}$ such that  $\forall f \in S_{n,j}, \quad
\exists f_i \geq f \in \bar{{\cal H}}_{n,j}$ satisfying
$ \tr{T_n(f)^{-1}T_n(f_i) - \id}/n \leq h_n(f_0,f_i)/8$ ,and
$ \tr{T_n(f_i-f)T_n^{-1}(f_0)}/ n \leq h_n(f_0,f_i)/8.$
Denote by $\bar{N}_{n,j}$ the logarithm
of the cardinality of the smallest possible net $\bar{{\cal H}}_n$.
$$ \bar{N}_{n,j} \leq n\varepsilon_n^2 j^\alpha, \quad \mbox{with}
\quad  \alpha <1.$$ Then, there exist $M,C,C'>0$ such that if
$\rho_n \leq \varepsilon_n^2$ and $n$ is large enough
\begin{eqnarray}\label{risk:rate}
E_0^n \left[ \pi \left. \left( f ; h_n(f_0,f) \geq M
\varepsilon_n^2 \right| X \right) \right] \leq \max \left(
e^{-n\varepsilon_n^2C}, \frac{ C^\prime }{ n ^2 } \right).
\end{eqnarray}

\end{thm}
\begin{proof}  
 Throughout the proof $C$ denotes a generic constant. We have
\begin{align}
 \pi& \left. \left( f ; h_n(f_0,f) \geq M
\varepsilon_n^2 \right| \Xn \right) = \frac{ \int\limits_{ f :
h_n(f_0,f) \geq M\varepsilon_n^2  }
\varphi_f(\Xn)/\varphi_{f_0}(\Xn) d\pi(f) }{
\int \varphi_f(\Xn)/\varphi_{f_0}(\Xn) d\pi(f)} \notag \\
=& \frac{ \int\limits_{ f : \varepsilon \geq h_n(f_0,f) \geq
M\varepsilon_n^2 } \varphi_f(\Xn)/\varphi_{f_0}(\Xn) d\pi(f) }{
\int \varphi_f(\Xn)/\varphi_{f_0}(\Xn) d\pi(f)}
+ \frac{ \int\limits_{ f :
h_n(f_0,f) \geq \varepsilon } \varphi_f(\Xn)/\varphi_{f_0}(\Xn)
d\pi(f) }
{ \int \varphi_f(\Xn)/\varphi_{f_0}(\Xn) d\pi(f)} \notag \\
=& \frac{ N_n }{ D_n } + R_{n,2},\notag
\end{align}
for some $\varepsilon >0$.
Theorem \ref{th:cons} implies that
$P_0\left[ R_{n,2}> e^{-n\delta} \right] \leq \frac{ C }{ n^2}$,
for some constants $C, \delta >0$. Then we consider the first term
of the right hand side of the above equation. Using an argument similar to the one used in the previous proof, let
 $$N_{n,j} = \int\limits_{ f :
\varepsilon_n^2j \leq h_n(f_0,f) \leq \varepsilon_n^2(j+1) }
\frac{ \varphi_f(\Xn) }{ \varphi_{f_0}(\Xn) }d\pi(f) $$
and
\begin{eqnarray*}
 E_0^n
\left[ \frac{ N_n }{D_n } \right] &\leq& \sum_{ j \geq M }
E_0^n\left[ \varphi_{n,j} \right] + E_0^n \left[
(1-\varphi_{n,j})\frac{ N_{n,j} }{D_n } \right],
\end{eqnarray*}
where $\varphi_{n,j} = \max_{f_i \in \bar{{\cal
H}}_{n,j}}\varphi_i$, and $\varphi_i$ is a test function defined as in the previous
Section, that is $\varphi_i = \I_{D_i}$, where

$$
D_i = \left \{ \Xn'(T_n^{-1}(f_i) - T_n^{-1}(f_0))\Xn
\geq \tr{\id - T_n(f_0)T_n^{-1}(f_i)} + h_n(f_0,f_i)/4\right \}.$$

 Then, (\ref{test:null:0}) implies that
 $$
 E_0^n\left[\phi_{n,j}\right] \leq  \sum_{i:f_i \in \bar{{\cal
 H}}_{n,j}} e^{-C n \varepsilon_n^2 j}
 \leq \bar{N}_{n,j} e^{-C n \varepsilon_n^2 j}
  \leq e^{-C n \varepsilon_n^2 j}.
$$
We also have that
\begin{eqnarray*}
E_0^n \left[ (1-\varphi_{n,j})\frac{ N_{n,j} }{D_n } \right] &\leq&
P_0^n\left[ D_n \leq e^{-n\rho_n}/2 \right] +
 2e^{n\rho_n}\pi(\bar{{\cal F}}_n^c\cap \{f:d(f,f_0) \leq \varepsilon\}) \\
   &+ &  2e^{n\rho_n} \int_{S_{n,j}}E_f^n\left[ 1 - \varphi_{n,j}
 \right]d\pi(f) \\
  &\leq&
2e^{-n\rho_n} + 2e^{n\rho_n}e^{-nC\varepsilon_n^2j^2} + P_0^n\left[
D_n \leq e^{-n\rho_n}/2 \right].
\end{eqnarray*}
Moreover, using the same calculations as in the proof of theorem
\ref{th:cons}
 \begin{eqnarray*}
 P_0^n\left[ D_n \leq e^{-n\rho_n}/2 \right] &\leq& P_0^n\left[
D_n \leq e^{-n\rho_n/2}\pi({\cal B}_n) \right] \\
 &\leq& \frac{\int_{\mathcal B_n} P_0^n\left[ \Omega_{n,1}^c(f) \right]d\pi(f)
 }{ \pi(\mathcal B_n)},
\end{eqnarray*}
where $\Omega_{n,1} = \{(\Xn,f); \Xn^t(T_n^{-1}(f)-T_n^{-1}(f_0))\Xn
- \log\det[A(f_0,f)] \leq n\rho_n/2 \}$. Using the exponential bound
similar to (\ref{test:null:0}), we obtain, if $f \in \mathcal B_n$, that
$$P_0^n\left[ \Omega_n^c \right] \leq \exp\{ -n \rho_n \left( \frac{\rho_n}{16
b_n(f_0,f)} \wedge \frac{ 1 }{ 8 } \right) \}
 \leq
e^{-n\rho_n/16},$$ on $\mathcal B_n$
 and Theorem \ref{th:rate:spde} is proved.
\end{proof}

The conditions given in Theorem \ref{th:rate:spde} are similar in spirit to
those considered for rates of convergence of the
posterior distribution in the i.i.d. case. The first one is a
condition on the prior mass of Kullback-Leibler neighborhoods of
the true spectral density, the second one is necessary to allow for
sets with infinite entropy (some kind of non compactness)  and the
third one is an entropy condition. The inequality
(\ref{risk:rate}) obtained in Theorem \ref{th:rate:spde} is non
asymptotic, in the sense that it is valid for all $n$. However,
the distances considered in Theorem \ref{th:rate:spde} heavily depend
 on $n$ and, although they express the impact of the
differences between $f$ and $f_0$ on the observations, they are
not of great practical use. For these reasons, the entropy condition is awkward and it cannot be directly transformed into some more common
entropy conditions. To state a result involving distances between
spectral densities that might be more useful, we need to consider some specific
class of priors, namely the FEXP priors, as defined in Section
\ref{sec:cons:fexp}. For this class we obtain rates of convergence in
terms of the $L_2$ distance between the logarithm of the spectral
densities, $l(f,f')$. The rates obtained are the optimal rates up
to a $\log n$ term, at least on certain classes of spectral
densities. It is to be noted that the calculations used when
working on these classes of priors are actually more involved
than those used to prove Theorem \ref{th:rate:spde}. This is
quite usual when dealing with rates of convergence of posterior
distributions, however this is emphasized here by the fact that
distances involved in Theorem 4 are strongly dependent on $n$.
The method used in the case of the FEXP prior can be extended to
other types of priors.

\subsection{The FEXP prior - rates of convergence}\label{sec:fexp:rate}

Here we apply Theorem \ref{th:rate:spde} to the FEXP
priors, which we define through a slightly different parameterization.
In particular,
$f(\lambda) =|1-e^{i\lambda}|^{-2d}\tilde{f}(\lambda)$, and
$\log{\tilde{f}}(\lambda)  = \sum_{j=0}^K \theta_j \cos{j\lambda}.$
Then the prior can be written in terms of a prior on $\left (
d,K,\theta_0,...,\theta_K \right )$.
Define now the classes of spectral densities
 $$ {\cal S}(\beta, L_0) = \{ h \geq 0; \log{h} \in L^2[-\pi,\pi], \log{h}(x) =
 \sum_{j=0}^\infty \theta_j \cos{j x}, \sum_j \theta_j^2 (1+j)^{2\beta} \leq L_0 \},
 $$
with $\beta >0$. Also, assume that there exists a real value $\beta > 0$
such that
$\tilde{f}_{0} \in {\cal L}^\star(M,m,L) \cap {\cal S}(\beta,
L_0)$. We can then write $f_0$ as
$$f_0(\lambda) = |1- e^{i\lambda}|^{-2d_0}\exp\left\{
\sum_{j=0}^\infty \theta_{j,0}\cos{j\lambda} \right \}.$$ Note that
$\beta$ is a smoothness parameter. Classes similar to ${\cal S}(\beta, L_0)$
are considered by
\cite{mousou:03}. We now describe the construction of the FEXP
prior, so that it can be adapted to ${\cal S}(\beta, L_0)$. Let
$S_K$ be a r.v. with density $g_A(\cdot)$, which is positive in the interval
$[0,A]$, let $\eta_j = \theta_j j^\beta$ and suppose that the prior
on $(\eta_1/S_K,...,\eta_K/S_K)$ has positive density on the set
$$\tilde{S}_{K+1} = \{x = (x_1,...,x_{K+1}); \sum_{j=1}^{K+1} x^2_j =
1\}$$.
We denote this class as the class of $\mbox{FEXP}(\beta)$
priors. Note that if $\beta > 1/2$ then there exists a constant $M$
and $\rho < 2\beta-1$ such that for all $f \in {\cal S}(\beta, L_0)$
associated with the parameters $(k,\theta_0,...,\theta_k)$ then
$$\sum_{i=0}^k|\theta_i| \leq M, \quad |\log{\tilde{f}}(x) -
\log{\tilde{f}}(y)| \leq M|x-y|^\rho.$$ first,
\begin{eqnarray*}
\sum_{i=0}^k|\theta_i| & \leq& \sum_{i=0}^ki^{2\beta} \theta_i^2  +
\sum_{i=0}^k|\theta_i| \I_{|\theta_i|\geq (1+i)^{2\beta}\theta_i^2} \\
 &\leq& L_0 + \sum_{i=0}^\infty (1+i)^{-2\beta} ,
 \end{eqnarray*}
and, second, since $\sum_j j^\rho |\theta_j|$
is uniformly bounded,
\begin{eqnarray*}
|\log{\tilde{f}}(x) - \log{\tilde{f}}(y)| &\leq& |x-y|^\rho
\left(2\sum_j|\theta_j|\right) \sum_{j\geq 1}|\theta_j| j^\rho \\
 &\leq& |x-y|^\rho \left(2\sum_j|\theta_j|\right) \sum_{j\geq 1}|\theta_j|
 j^\rho \\
  &\leq& M|x-y|^\rho.
\end{eqnarray*}
Therefore the prior lies in $\bar{\mathcal G}(t,m,M,L,\rho)$ for
some positive constant $m,M,L,\rho$.

 We now give the rates of convergence associated with the
$\mbox{FEXP}(\beta)$ priors,
 when the true spectral density belongs to $ {\cal S}(\beta,
 L_0)$.
\begin{thm}\label{th:fexp:rate}
Assume that there exists $\beta >\frac12$ s.t.
$\tilde{f}_{0} \in {\cal L}^\star(e^{L_0},e^{-L_0},L) \cap {\cal S}(\beta,
L_0)$.  Let $\pi$ be a $\mbox{FEXP}(\beta)$ prior and assume that
\bei
\item[i)]
$K\sim\mbox{ Poi}(\mu)$;
\item[ii)] the prior on $d$ is strictly positive on $[-1/2+t,1/2-t]$,
with $t>0$;
\item[iii)] $S_K$ has a positive density on $(0,A)$
with $A$
such that $A^2 \geq L_0$.
\eni
Then there exist $C, C'>0$ such that, for $n$ large enough
\begin{eqnarray} \label{post:rate}
P^\pi\left[ \{ f \in {\cal F}^+ : l(f,f_0) > C n^{-2\beta/(2
 \beta+1)} \log{n}^{(2\beta+3)/(2\beta+1)} \}| \Xn \right] \leq \frac{ C' }{ n^2 }
 \end{eqnarray}
and
\begin{eqnarray}\label{risk:rate2}
 E_0^n\left[l(\hat{f},f_0) \right] \leq  2C  n^{-2\beta/(2
 \beta+1)} \log{n}^{(2\beta+3)/(2\beta+1)},
 \end{eqnarray}
where $\log{\hat{f}} (\lambda)  = E^\pi \left[ \log{f}(\lambda)|\Xn \right]$.
\end{thm}

\begin{proof}
Throughout the proof, $C$ denotes a generic constant. The proof of the
theorem is divided in two parts; in the first
part, we prove that
\begin{eqnarray}\label{fexp:rate:1}
E_0^n \left[ P^\pi\left\{ f : h_n(f,f_0) \geq
n^{-2\beta/(2\beta+1)} \log{n}^{(2\beta+3)/(2\beta+1)} | \Xn
\right\} \right] \leq \frac{ C }{ n^2 },
\end{eqnarray}
while in the second part we prove that
\begin{eqnarray}\label{fexp:rate:2}
h_n(f,f_0) \leq Cn^{-\frac{2\beta}{2\beta+1}} \log{n}^{1/\beta}
\Rightarrow l(f,f_0) \leq C'n^{-\frac{2\beta}{2\beta+1}}
\log{n}^{\frac{2\beta+3}{2\beta+1}},
\end{eqnarray}
 for some constant $C'>0$, when $n$ is large enough.
 The latter inequality implies that
\begin{eqnarray*}
E_\pi \left[ l(f,f_0) |\Xn \right] &\leq &
C'n^{-\frac{2\beta}{2\beta+1}} \log{n}^{\frac{2\beta+3}{2\beta+1}} +
\int\limits_{A(n,\beta)} l(f,f_0)d\pi(f|\Xn) \\
 &\leq&
 2C'n^{-\frac{2\beta}{2\beta+1}}\log{n}^{\frac{2\beta+3}{2\beta+1}},
\end{eqnarray*}
for large $n$, where
$A(n,\beta)=\{h_n(f,f_0)>
Cn^{-\frac{2\beta}{2\beta+1}}\log{n}^{\frac{2\beta+3}{2\beta+1}}\}$. This would imply
Theorem \ref{th:fexp:rate}.
To prove (\ref{fexp:rate:1}), we need to show that conditions
(i)-(iii) of Theorem \ref{th:rate:spde} are fulfilled. Condition
(ii) is obvious because the prior has the same form as in Section
\ref{sec:cons:fexp} and, because $\mathcal{S}(\beta,L)  \subset
\bar{\mathcal G}(t,m,M,L',\rho)$, with $t,m,L',\rho$ positive
constant depending on $\beta,L$. Thus we can choose
 $$\bar{\mathcal F}_n = \left \{ f(\lambda) = |1-e^{i\lambda}|^{-2d}\exp\left (
 \sum_{j=0}^K \theta_j \cos{(j\lambda)}\right ) \right \},$$
with $ K \leq K_n, |d-d_0| \leq \delta,
  \sum_jj^{2\beta}\theta_j^2 \leq L_0$, leading to
$$ \pi\left(\bar{\mathcal F}_n \cap \{ f,h(f,f_0)< \epsilon\}
\right) \leq \pi(K \geq K_n ) \leq e^{-K_n\log{K_n}}$$ for $K_n$
large enough. By choosing $K_n = k_0
n^{1/(2\beta+1)}\log{n}^{2/(2\beta+1)}$, we obtain
$$ \pi\left(\bar{\mathcal F}_n \cap \{ f,h(f,f_0)< \epsilon\}
\right) \leq  e^{-k_0
n^{1/(2\beta+1)}\log{n}^{(1-2\beta)/(2\beta+1)}}.$$ Hence, letting
$\rho_n  = \epsilon_n^2 = n^{-2\beta/(2\beta +1)}
\log{n}^{(2\beta+3)/(2\beta+1)}$, condition (ii) is satisfied.
We now show that assumption (i) of Theorem \ref{th:rate:spde} is
satisfied. Let  $d \leq d_0 \leq d+ \epsilon_n/\log{n}^{3/2}$ and,
for all $l=0,\dots,K_n$,
$$|\theta_l - \theta_{0l}| \leq  (l+1)^{-(\beta+1/2)}(\log{(l+1)})^{-1}\epsilon_n/\log{n}^{3/2}.$$
 Since
 $f_0 \in \mathcal{S}(\beta,L_0)$, $\exists\, t_0>0$ such that
 \begin{eqnarray}\label{eq:f0n}
 \sum_{l \geq K_n} \theta_{0l}^2 \leq L_0 K_n^{-2\beta} \leq C \epsilon_n^2
 (\log{n})^{-3}, \quad \sum_{l \geq K_n} |\theta_{0l}| \leq
 K_n^{-t_0}
\end{eqnarray}
 Since
\begin{eqnarray*}
 KL_n(f_0;f) &\leq& h_n(f_0,f) \\
 &=&  \frac{ 1 }{ 2n } \tr{T_n(f_0-f)T_n^{-1}(f)T_n(f_0-f) T_n^{-1}(f_0)},
\end{eqnarray*}
it is enough to prove the assumption under the above conditions
for $h_n(f,f_0) \leq C\epsilon_n^2$. The difficulty here comes from
the strong dependence on $n$ of the distance $h_n$.
Let
$$f_{0n}(\lambda) = |1 - e^{i\lambda}|^{-2d_0}\exp{\left(\sum_{l=0}^{K_n}\theta_{0l}\cos{l\lambda}\right)},
\quad  b_n(\lambda) = 1 - \exp\left(-\sum_{l \geq K_n+1}
\theta_{l0}\cos{l\lambda}\right),$$
and
$g_n =f_{0n}^{-1}(f_{0n}-f)$;
then $f_0- f = f_0b_n + f_{0n}g_n$ and
\begin{eqnarray}\label{upbound:hn1}
 nh_n(f_0,f)\!\! &\leq& \!\!
 \tr{T_n(f_0b_n)T_n^{-1}(f)T_n(f_0b_n) T_n^{-1}(f_0)
}  \nonumber \\
  & &  +  \tr{T_n(f_{0n}g_n)T_n^{-1}(f)T_n(f_{0n}g_n) T_n^{-1}(f_0)}.
\end{eqnarray}
Both terms of the right hand side of  (\ref{upbound:hn1}) are
treated similarly, using Lemma \ref{whittles} we can bound them by
 \begin{eqnarray*}
 \tr{ T_n(f_{0n}b_n)T_n^{-1}(f)T_n(f_{0n}b_n)
T_n^{-1}(f_0)}
 &\leq&   C(\log{n})^3  n|b_n|_2^2 + O(n^{\delta}). \\
  \tr{ T_n(f_{0n}b_n)T_n^{-1}(f)T_n(f_{0n}b_n)
T_n^{-1}(f_0)}
 &\leq&   C(\log{n})^3  n|g_n|_2^2 + O(n^{\delta}).
\end{eqnarray*}
This implies that $h_n(f_0,f) \leq C \epsilon_n^2$, when $f$
satisfies the conditions described above and
\begin{eqnarray*}
\mathcal{B}_n &\subset& \left\{ f_{k,d,\theta}; k \geq
 K_n, d \leq d_0 \leq d + \frac{ \epsilon_n }{ (\log{n})^{3/2} }, 0 \leq l \leq
 K_n,\right .\\
 &&\left . |\theta_l- \theta_{0l}| \leq \frac{ (l+1)^{-(\beta+1/2)}\epsilon_n }{ (\log{(l+1)})\log{n}^{3/2} } \right\}.
 \end{eqnarray*}
 The prior
probability of the above set is bounded from below by
 $$\pi(K_n) \mu_1 \left( (\eta_1,...,\eta_{K_n}):
 |\eta_l-\eta_{0l}| \leq C\frac{ l^{-1/2}\epsilon_n }{ (\log{l})\log{n}^{3/2} } \right)\rho_n \log{n}^{-3/2},$$
where $\mu_1$ denotes the uniform measure on the set $\{
(\eta_1,...,\eta_{K_n}); \sum_l \eta_l^2 \leq A \}$. We finally
obtain that
 $$\pi(\mathcal{B}_n(\delta)) \geq e^{-CK_n \log{n}} \geq
 e^{-n\rho_n/2}$$ by choosing $k_0$ small enough, and condition
 (i) of Theorem \ref{th:fexp:rate} is satisfied by the
 FEXP($\beta$) prior. We now verify condition (iii) of Theorem
 \ref{th:fexp:rate}.
  Let $j_0 \leq j \leq J_n$, where $j_0$ is some positive constant, and consider $f \in S_{n,j}$,
  as defined in Theorem \ref{th:rate:spde},    where
 $f(\lambda) = f_{\theta,k} = |1 - e^{i\lambda}|^{-2d}\exp\{ \sum_{l=1}^k
 \theta_l \cos{(l \lambda)} \}$.
Consider
  $$f_u(\lambda) = |1 - e^{i\lambda}|^{-2d_u }\exp\{
  \sum_{l=1}^k \theta_l^u \cos{(l \lambda)} \};  f(\lambda)=|1 -
e^{i\lambda}|^{-2d }\exp\{ \sum_{l=1}^k \theta_l \cos{(l \lambda)}
\}, $$ then, if $c\geq c_0\sum_{l\geq 0}
(l+1)^{-\beta-1/2}\log{(l+1)}^{-1}$, $c_0>0$ and
 $d^u \geq d \geq d^u - c\epsilon_n^2 j$, \, $|\theta_l -\theta_l^u| \leq
 c_0(l+1)^{-\beta-1/2}\log{(l+1)}^{-1}\epsilon_n^2 j,$ \, $l \geq 1,$ \,
$\theta_0^u - 4c\epsilon^2_n j \leq \theta_0 \leq \theta_0^u -
 3c\epsilon_n^2j$, one obtains
$$1 \leq  \frac{ f_u }{ f }(\lambda) \leq (1 - \cos{\lambda})^{-2c\epsilon_n^2j}e^{5 c \epsilon_n^2
j}$$
  and
\begin{eqnarray*}
\tr{ T_n^{-1}(f)T_n(f_u-f)} &\leq& 15 c \epsilon_n^2j
 \tr{T_n^{-1}(f) T_n(f_u)} \\
 &\leq& Cc \epsilon_n^2 j \\
  &\leq & Cc h_n(f_0,f_u).
\end{eqnarray*}
Choosing $c$ small enough one obtains that
$\tr{T_n^{-1}(f)T_n(f_u-f)} \leq n h_n(f_0,f_u)/8$. Similarly
\begin{eqnarray*}
\tr{T_n^{-1}(f_0)T_n(f_u-f)} &\leq& 4 c
\epsilon_n^2j  \tr{T_n^{-1}(f_0)T_n(f_u)} \\
 &\leq&  cC h_n(f_0,f_u)/8.
\end{eqnarray*}
Since we are in the set $\{ f; h(f_0,f) \leq \epsilon \}$, for some
$\epsilon >0$ fixed but as small as we need, there exists $\epsilon'
, \epsilon">0$ such that
 $$|d-d_0| < \epsilon', \quad \sum_{l=1}^K(\theta_l-
\theta_{l0})^2 + \sum_{l \geq K+1}\theta_{l0}^2 \leq \varepsilon".$$
Let $K \leq K_n = K_0 n^{1/(2\beta+1)}(\log{n})^{-1}$, the number of
$f_u$ defined as above in the set $S_{n,j}$ is bounded by
 $$ N_{n,j} \leq K_n j^{-1} \epsilon_n^{-2}\left( C K_n j^{-1} \epsilon_n^{-2} \right)^{K_n}
 $$
and
$$\bar{N}_{n,j} = \log{N_{n,j}} \leq c_j n\epsilon_n^2$$
 where $c_j$ is decreasing in $j$. Hence by choosing $j_0$
large enough condition (iii) is verified by the FEXP($\beta$)
prior. This achieves the proof of (\ref{fexp:rate:1}) and we
obtain a rate of convergence, in terms of the distance $h_n(\cdot,\cdot)$.
We now prove (\ref{fexp:rate:2}) to obtain a rate of convergence
in terms of the distance $l(\cdot,\cdot)$. Consider $f$ such that
 \begin{eqnarray*}
  h_n(f_0,f) &=&
 \frac{ 1 }{ n }\tr{T_n^{-1}(f_0)T_n(f-f_0)T_n^{-1}(f)T_n(f-f_0)} \leq \epsilon_n^2.
 \end{eqnarray*}
 Equation (\ref{whit01}) of Lemma \ref{whittles} implies that
  \begin{eqnarray*}
 \frac{ 1 }{ n }\tr{ T_n(f_0^{-1})T_n(f-f_0)T_n(f^{-1})T_n(f-f_0)}
 \leq  C\epsilon_n^2,
 \end{eqnarray*}
leading to
\begin{eqnarray}\label{res:1}
\frac{ 1 }{ n }\tr{T_n(g_0)T_n(f-f_0)T_n(g)T_n(f-f_0)} \leq
C\epsilon_n^2,
 \end{eqnarray}
where $g_0 = (1-\cos{\lambda})^{d_0}$, $g = (1-\cos{\lambda})^{d}$.\\
We now prove that  $\tr{T_n(g_0(f-f_0))T_n(g(f-f_0))} \leq
C\epsilon_n^2$: we use the same representation as in the treatment
of $\gamma(b)$ in Appendix \ref{app:lemmas}. For the sake of simplicity
we consider the case $d \geq d_0$
\begin{eqnarray*}
\bar{\Delta}  &= & \frac{ 1 }{ n } \tr{
T_n(g_0(f-f_0))T_n(g(f-f_0))} \\
&-& \frac{ 1 }{n }
\tr{T_n(g_0)T_n(f-f_0)T_n(g)T_n(f-f_0)} \\
 &=&
 \frac{ 1 }{ n } \int\limits_{[-\pi,\pi]^3}
(f-f_0)(\lambda_2)g_0(\lambda_2) (f-f_0)(\lambda_4)g(\lambda_4)
\left( \frac{ g_0(\lambda_1) }{ g_0(\lambda_2) } - 1 \right)\\
&\times& \Delta_n(\lambda_1-\lambda_2)\Delta_n(\lambda_2-\lambda_4)
 \Delta_n(\lambda_4-\lambda_1) d\underline{\lambda} \\
 &+ &
 \frac{ 1 }{ n } \int\limits_{[-\pi,\pi]^4}
(f-f_0)(\lambda_2)g_0(\lambda_1) (f-f_0)(\lambda_4)g(\lambda_4)
\left( \frac{ g(\lambda_3) }{ g(\lambda_4) } - 1 \right)\\
  & \times&   \Delta_n(\lambda_1-\lambda_2)\Delta_n(\lambda_2-\lambda_3)
 \Delta_n(\lambda_3-\lambda_4)\Delta_n(\lambda_4-\lambda_1)
d\underline{\lambda} \\
\end{eqnarray*}
\begin{eqnarray*}
 &\leq&  \frac{ C \log{n} }{ n } \int\limits_{[-\pi,\pi]^2}
 |\lambda_2|^{-2(d-d_0)}|\lambda_1|^{-1+\delta}
L_n(\lambda_1-\lambda_2)^\delta L_n(\lambda_2-\lambda_1)
d\underline{\lambda} \\
 &+ &
 \frac{ C }{ n } \int\limits_{[-\pi,\pi]^4} \frac{|\lambda_1|^{2d}|}{|\lambda_2|^{2d}
 \lambda_3|^{1-\delta}}
L_n(\lambda_1-\lambda_2)L_n(\lambda_2-\lambda_3)L_n(\lambda_3-\lambda_4)^\delta
L_n(\lambda_4-\lambda_1) d\underline{\lambda} \\
 &\leq&
 \frac{ C \log{n}n^{2\delta} }{ n }\\
&+&
 \frac{ C\log{n}n^\delta }{ n } \int\limits_{[-\pi,\pi]^3}
\frac{|\lambda_1|^{2d}|}{|\lambda_2|^{2d} \lambda_3|^{1-\delta}}L_n(\lambda_1-\lambda_2)L_n(\lambda_2-\lambda_3)L_n(\lambda_3-\lambda_1)
d\underline{\lambda} \\
 & \leq& \frac{ C (\log{n})^{2} }{ n^{1 - 2\delta}},
 \end{eqnarray*}
if $\delta \geq 4(d-d_0)$. We have used inequality (\ref{ineq:Dahl})
together with inequality (\ref{ineq:int:Ln}).
 This implies, together with
(\ref{res:1}) that
 $$\frac{ 1 }{ n } \tr{ T_n(g_0(f-f_0))T_n(g(f-f_0))}  \leq C \epsilon_n^2. $$
To finally obtain (\ref{fexp:rate:2}), we use equation
(\ref{whit02}) in Lemma \ref{whittles} which implies that
\begin{eqnarray*}
A_n &=& \tr{ T_n(g_0(f-f_0))T_n(g(f-f_0))} - \tr{
T_n(g_0g(f-f_0)^2)} \\
 & \leq & C n^{-1+\delta} + \log{n} \sum_{l=0}^{K_n} l|\theta_l|
\left(
  \int_{[-\pi,\pi]}g_0g(f-f_0)^2(\lambda)d\lambda\right)^{1/2}.
\end{eqnarray*}
Moreover
 \begin{eqnarray*}
 \sum_{l =1}^{K_n} l |\theta_l| &\leq & \sum_{l=1} l^{2\beta+r}
 \theta^2 + \sum_{l=1}^{K_n}l^{-r/(2\beta-1)} \\
  &\leq& C K_n^r + K_n^{1-r/(2\beta-1)},
 \end{eqnarray*}
by choosing $r= (2\beta-1)/2\beta$, $A_n/n$  is of order
$n^{-(4\beta^2 +1)/(2\beta(2\beta+1))}$ which is negligible
compared to $n^{-2\beta/(2\beta+2)}$ so that if $ \beta \geq 1/2$
 $$\int_{[-\pi,\pi]} g_0g (f_0-f)^2 d\lambda \leq \epsilon_n^2,$$
 which achieves the proof.
\end{proof}

\section{Discussion}
In this paper we have considered the theoretical properties of our Bayesian procedure. A related and important problem, which deserves the same attention,
is its practical implementation.
Due to the length of the present paper, we discuss these issues elsewhere; see for example \cite{lr:06}; here we only sketch the main features of the proposed algorithm.

From a computational perspective, the practical implementation of a nonparametric Bayesian analysis based on a FEXP prior and a Gaussian likelihood, is plagued by two difficulties: i) the number of parameters to estimate varies with $K$ the number of terms in the FEXP expansion; ii) the likelihood function is quite expensive to evaluate, due to the Toeplitz structure of the covariance matrix.

After trying several approaches we finally recommend the use of the $D$-kernel Population MonteCarlo algorithm, presented and discussed in \cite{douc:05}, and which can be easily adapted to the varying dimension set-up.
For the evaluation of the inverse and of the determinant of the Toeplitz covariance matrix, we have used the algorithms proposed in
\cite{Hur:2006}.

\appendix
\section{\texorpdfstring{Lemmas \ref{lem:unifTR1} and \ref{lem:unifTR2}}{Lemmas 1 and 2}}
We state two technical lemmas, which are extensions of
\cite{jud:03} on uniform convergence of traces of Toeplitz
matrices, and which are repeatedly used in the paper.
\begin{lemma} \label{lem:unifTR1}
Let $t>0$, $M >0$  and $\bar{M}$ a positive function on $]0,\pi[$, let $p$ be a positive integer,
and
 $$\tilde{\tilde{\mathcal{F}}}(d,M,\bar{M}) = \left\{ f \in \tilde{F}, \forall u >0,  \sup_{|\lambda| > u}
  \frac{ d\tilde{f}(\lambda) }{ d\lambda } \leq \bar{M}(u) \right\},$$
  we have:
\begin{eqnarray}\label{unif1}
\sup_{\substack{p(d_1+d_2) \leq 1/2-t \\ f_i \in \tilde{\tilde{\mathcal{F}}}(d_1,M,\bar{M}) \\
g_i \in \tilde{\tilde{\mathcal{F}}}(d_2,M,\bar{M})}}
\left|
 \frac{ 1 }{ n }
 \tr{\prod_{i=1}^pT_n(f_i)T_n(g_i) }-
 (2\pi)^{2p-1}\int\limits_{-\pi}^\pi \prod_{i=1}^p
 f_i(\lambda)g_i(\lambda) d \lambda \right| \rightarrow 0.
\end{eqnarray}
and let $L>0 $ and $\rho \in (0,1]$
\begin{eqnarray}\label{unif2}
\sup_{\substack{p(d_1+d_2) \leq 1/2-t \\ f_i \in \mathcal{F}(d_1,M,L, \rho) \\
g_i \in \mathcal{F}(d_2,M,L, \rho)}}
\left|
 \frac{ 1 }{ n }\tr{\prod_{i=1}^pT_n(f_i)T_n(g_i) } -
 (2\pi)^{2p-1}\int\limits_{-\pi}^\pi \prod_{i=1}^p
 f_i(\lambda)g_i(\lambda) d \lambda \right| \rightarrow 0.
\end{eqnarray}
\end{lemma}

\noindent
This lemma is an obvious adaptation from \cite{jud:03}, and the only non obvious part is the change from the
condition of continuous differentiability in that paper to the
Lipschitz condition of order $\rho$, considered equation \ref{unif2}.
This different assumption affects only  equation (30) of \cite{jud:03}, with
$\eta_n$  replaced by $\eta_n^\rho$, which does not change the convergence results.

\begin{lemma}\label{lem:unifTR2}
\begin{eqnarray*}
\sup_{\substack{2p(d_1-d_2) \leq \rho_2\wedge 1/2 -t \\ f_i \in {\cal F}(d_1,M,L,\rho_1) \\
g_i \in {\cal G}(d_2,m,M,L,\rho_2)}} \left|
 \frac{ 1 }{ n }\tr{\prod_{i=1}^pT_n(f_i)T_n(g_i)^{-1} } -
 \frac{ 1 }{ 2\pi }\int_{-\pi}^\pi \prod_{i=1}^p\frac{ f_i(\lambda) }{ g_i(\lambda) } d \lambda \right| \rightarrow
 0,
\end{eqnarray*}

\begin{eqnarray*}
\sup_{\substack{2p(d_1-d_2) \leq \rho_2\wedge 1/2 -t \\ f_i \in \tilde{\tilde{\mathcal{F}}}(d_1,M,\bar{M}) \\
g_i \in {\cal G}(d_2,m,M,L,\rho_2)}} \left|
 \frac{ 1 }{ n }\tr{\prod_{i=1}^p T_n(f_i)T_n(g_i)^{-1} } -
 \frac{ 1 }{ 2\pi }\int_{-\pi}^\pi \prod_{i=1}^p\frac{ f_i(\lambda) }{ g_i(\lambda) } d \lambda \right| \rightarrow
 0.
\end{eqnarray*}
and
\begin{eqnarray*}
\sup_{\substack{2p(d_1-d_2) \leq 1/2-t \\ f_i \in \tilde{\tilde{\mathcal{F}}}(d_1,M,\bar{M}) \\
g_i \in {\cal L}(d_2,m,M,L)}} \left|
 \frac{ 1 }{ n }\tr{\prod_{i=1}^pT_n(f_i)T_n(g_i)^{-1} } -
 \frac{ 1 }{ 2\pi }\int_{-\pi}^\pi \prod_{i=1}^p\frac{ f_i(\lambda) }{ g_i(\lambda) } d \lambda \right| \rightarrow
 0.
\end{eqnarray*}

\end{lemma}

\begin{proof}
In this second lemma, the uniformity result is a consequence of
the first lemma, as in \cite{jud:03};
 The only difference is in the proof of Lemma 5.2. of \cite{dhl:89}, i.e. in the study of terms in the form
 $$ |\id -T_n(g)^{1/2}T_n\left((4\pi^2 g)^{-1}\right)T_n(g)^{1/2}|.$$
  Following Dahlhaus's \cite{dhl:89} proof, we obtain an upper bound of  $$\left| \frac{ g(\lambda_1) }{ g(\lambda_2) } - 1 \right|$$
  which is different  from \cite{dhl:89}.
  If $ g \in {\cal G}(d_2,m,M,L,\rho_2)$, the Lipschitz condition in $\rho$ implies
   that
 $$
 \left| \frac{ g(x) }{ g(y) } - 1 \right| \leq
 K\left( |x-y|^\rho + \frac{ |x-y|^{1-\delta} }{ |x|^{1-\delta} }
 \right).$$
 Calculations using $L_N$ as in \cite{dhl:89} imply that
$$ |I - T_n(f)^{1/2}T_n\left((4\pi^2
 f)^{-1}\right)T_n(f)^{1/2} |^2 =O(n^{1-2\rho}\log{n}) +
 O(n^\delta), \quad \forall \delta >0. $$
If $g \in \mathcal{L}^\star(M,m,L)$ as defined in Section
\ref{sec:fexp:rate}, then
$$
 \left| \frac{ f(x) }{ f(y) } - 1 \right| \leq
 K\left(  \frac{ |x-y|^{1-3\delta} }{ (|x| \wedge |y|) ^{1-\delta}
 } \right) \leq K  |x-y|^{1-3\delta}\left( \frac{ 1 }{ |x|^{1-\delta}} + \frac{ 1 }{ |y|^{1-\delta} } \right)$$
and \cite{dhl:89} Lemma 5.2 is proved, leading to a
constraint in the form $4p(d_1-d_2)< 1$ (corresponding to
$\rho=1$).

 \noindent
Then, using again Dahlhaus' (1989) calculations, we
obtain that
 $$|A - B| =0(n^{2(d_2-d_1)}n^{1/2 - (\rho \wedge 1/2)+\delta}), \quad \forall \delta > 0$$
 and finally that
\begin{eqnarray*}
\frac{ 1 }{ n } \tr{ \prod_{j=1}^p A_j - \prod_{j=1}^p B_j}
  & =& \sum_{k=1}^p O(n^{-1/2}n^{2(p-k)(d_2 - d_1)}n^{2(d_2-d_1)}n^{1/2 - \rho}) \\
  & = &\sum_{k=1}^p O(n^{2(p-k+1)(d_2-d_1)-(\rho\wedge 1/2) })
\end{eqnarray*}
   which goes to $0$ when $2p(d_2-d_1)< \rho \wedge 1/2$.
   \end{proof}

\section{\texorpdfstring{Proof of Theorem \ref{th:cons}}{Proof of Theorem 2.1}}
\label{app:pr:th:cons}
Before giving the proof of Theorem \ref{th:cons}, we give a few
notations that are used throughout the paper: Let $f,f_1$ be
spectral densities:
  \bei
  \item
  $A(f_1,f) = T_n(f)^{-1}T_n(f_1)$
  \item
  $B(f_1,f) =T_n(f_1)^{1/2}[T_n(f)^{-1} - T_n(f_1)^{-1}]T_n(f_1)^{1/2}$
  \item
  $b_n(f_1,f) = \tr{\id - T_n(f_1)T_n(f)^{-1})^2}/n$
  \item
  $b(f_1,f) =  (2\pi)^{-1} \int_{-\pi}^\pi (f_1/f - 1)^2(x)dx$.
 \eni

 {\it Proof of Theorem \ref{th:cons}}

 The proof follows the same ideas as in \cite{ggvdv:01}. The
 main difficulty here is to transform constraints on quantities such
 as $h_n(f,f_0)$ or $KL_n(f,f_0)$ in terms of distances between
 $f,f_0$ independent on $n$, uniformly over $f$.

 We can write
 \begin{eqnarray*}
 P^\pi\left[ A_\varepsilon^c |\Xn \right] &=& \frac{
 \int_{A_\varepsilon^c} \varphi_f(\Xn)/\varphi_{f_0}(\Xn)d\pi(f) }{  \int_{\tilde{{\cal F}}_+}
 \varphi_f(X)/\varphi_{f_0}(\Xn)d\pi(f) } =   \frac{ N_n }{ D_n }.
\end{eqnarray*}
Then the idea is to bound from below the denominator using
condition (i) of the Theorem and to bound from above the numerator
using a discretization of $A_\epsilon$ based on the net ${\cal
H}_n$ defined in (ii) of the Theorem and on tests.

Let $\varepsilon > \delta >0$: one has
\begin{eqnarray}\label{pr:thcons1}
P_0\left[ P^\pi\left[ A_\varepsilon^c |\Xn \right] \geq e^{-n\delta
} \right] &\leq& P_0^n\left[ D_n \leq e^{-n\delta} \right] +
P_0^n\left[ N_n \geq e^{-n 2\delta} \right] \nonumber \\
 &=& p_1 + p_2
\end{eqnarray}
Also, let
 $$\tilde{{\cal B}}_n(c) = \{ f \in \bar{{\cal
G}}(t,M,m,L,\rho): nKL_n(f_0,f) \leq n c  \}, \quad  c>0 .$$ Using
Lemma 2, when $n$ is large enough,
$$\tilde{{\cal B}}_n(\delta/2) \supset \{ f \in \bar{{\cal
G}}(t,M,m,L,\rho); h(f_0;f) \leq \frac{\delta}{4}, 8(d_0 - d) \leq
\rho \wedge 1/2 - t \} = \mathcal B_{\delta/4}$$
 so that assumption (i) implies that, for $n$ large enough,
$$\pi(\tilde{{\cal B}}_n(\delta/2)) \geq \pi(\mathcal B_{\delta/4})
\geq e^{-n \delta/2}/2.$$
Define
$$\Omega_n = \{ (f,X) :  -X^t[T_n(f)^{-1}
  - T_n(f_0)^{-1}]X + \log{(\det(A(f_0,f)))} > -n\delta \}.$$
We then have
\begin{eqnarray*}
p_1 &\leq&  P_0^n\left( \int_{\Omega_n\cap \tilde{{\cal B}}_n}
\frac{ \varphi_f(X) }{\varphi_{f_0}(X) }d\pi(f) \leq
e^{-n\delta/2}\frac{ \pi(\tilde{{\cal B}}_n) }{ 2 } \right) \\
 &\leq&
P_0^n\left( \pi(\tilde{{\cal B}}_n \cap \Omega_n ) \leq \frac{
\pi(\tilde{{\cal B}}_n) }{2 } \right) \\
 &\leq&
P_0^n\left( \pi(\tilde{{\cal B}}_n \cap \Omega_n^c )> \frac{
\pi(\tilde{{\cal B}}_n)  }{2 } \right) \\
 &\leq&
2 \frac{ \int_{\tilde{{\cal B}}_n}P_0^n[\Omega_n^c ] d\pi(f) }{
\pi(\tilde{{\cal B}}_n) }.
\end{eqnarray*}
Moreover,
\begin{eqnarray*}
P_0^n[\Omega_n^c ] &=& P_0^n\left( \Xn^t[T_n(f)^{-1}
 - T_n(f_0)^{-1}]\Xn - \log{(\det(A(f_0,f)))} > n\delta \right) \\
  &=&
  Pr[y^tB(f_0,f)y - \tr{B(f_0,f)}\\
  &>& n\delta + \log{(\det(A(f_0,f)))} - \tr{B(f_0,f)}],
\end{eqnarray*}
where $y\sim N_n(\bzero, \id)$. \noindent When $f \in \tilde{{\cal
B}}_n$,
$n\delta + \log{(\det(A(f_0,f)))} - \tr{B(f_0,f)} > n\delta/2$,\\
so that
\begin{eqnarray*}
P_0^n[\Omega_n^c ] &\leq&  Pr[y^tB(f_0,f)y - \tr{B(f_0,f)} >
n\delta/2 ] \\
 &\leq&
 \frac{ 4E[(y^tB(f_0,f)y - \tr{B(f_0,f)})^4 ] }{ n^4 \delta^4 } \\
  &\leq& \frac{ \tr{B(f_0,f)^4} C }{ n^3 \delta^4}.
\end{eqnarray*}
Therefore, for all $f \in \tilde{{\cal B}}_n$,
\begin{eqnarray*}
P_0^n[\Omega_n^c ] &\leq& \frac{ M'C }{ n^3 \delta^4},
\end{eqnarray*}
and \begin{eqnarray} \label{p1:fin} p_1 \leq C_1/n^3 ,
\end{eqnarray} where $C_1$ is a positive constant.

We now consider the second term of  (\ref{pr:thcons1}), namely:
\begin{eqnarray*}
p_2 &=& P_0^n\left[ N_n \geq e^{-2n\delta} \right] \\
 & \leq& 2 e^{2n\delta} \pi({\cal F}_n^c ) + P_0^n\left[
\int_{A_\varepsilon^c \cap {\cal F}_n} \frac{ \varphi_f(\Xn) }{
\varphi_{f_0}(\Xn) } d\pi(f) \geq e^{- 2n\delta }/2 \right] \\
 &\leq & e^{-n(r-2\delta)} + \tilde{p}_2,
\end{eqnarray*}
take $2\delta <r$ and consider $\tilde{p}_2$. Consider the following
tests : let $f_i \in {\cal H}_n$,
 $$\phi_i = \I_{ X'(T_n^{-1}(f_0) - T_n^{-1}(f_i) )X \geq n\rho_i
 }.$$

 Recall that $\gamma = \rho_0 \wedge \rho \wedge 1/2 - t$ (or $ \rho \wedge 1/2 -
 t$, $\rho_0  \wedge 1/2 - t$, $1/2 - t$ depending on whether the
 spectral densities belong to ${\mathcal G}$ or $\bar{ \mathcal
 L}$).
We now prove that $E_0^n[\phi_i] \leq e^{-n\varepsilon
|\log{\varepsilon}|^{-1}}$ and $E_f^n[1-\phi_i] \leq
e^{-n\varepsilon |\log{\varepsilon}|^{-1}}$ for $f$ close to $f_i$.
\begin{enumerate}
 \item  If $4|d_0-d_i| \leq \gamma$, set
 $\rho_i = \tr{\id - T_n(f_0)T_n^{-1}(f_i)}/n + h_n(f_0,f_i)$,
then for all $1/4 > s>0$,
 \begin{eqnarray*}
&E_0^n[\phi_i]& \leq \exp\left\{-s n\rho_i\right\}
 E_0^n\left[ e^{s\Xn'(T_n^{-1}(f_0)- T_n^{-1}(f_i) )\Xn} \right] \\
  &=& \exp\left\{-s n\rho_i\right \}
  \exp\{- \frac{ 1 }{ 2 } \log{\det[ \id +
  2sB(f_0,f_i) ]}\} \\
  &\leq& \exp\left\{ -s n \rho_i -s \tr{ B(f_0,f_i)} +  s^2 \tr{ ((\id +
  s\tau B(f_0,f_i))^{-1}B(f_0,f_i))^2} \right \} \\
  &\leq& \exp\left\{ -s n \rho_i -s \tr{ B(f_0,f_i)} + 4 s^2
  \tr{B(f_0,f_i)^2}\right \},
\end{eqnarray*}
where $\tau \in (0,1)$ and the latter inequality is due to
 $$\id +   s\tau B(f_0,f_i) = \id (1-2s\tau) + 2s\tau A(f_0,f_i)
 \geq \frac{ 1 }{ 2 }\id , \quad \mbox{if } s<1/4.$$
 Replacing $\rho_i$ by its above expression and choosing $s$ to optimize the latter expression, we obtain:
\begin{eqnarray}   \label{test:null:0}
E_0^n[\phi_i] \leq   \max\left( \exp{\{-n\frac{ h_n(f_0,f_i)^2 }{
16 b_n(f_0,f_i) } \}},\exp{\{-n\frac{ h_n(f_0,f_i) }{ 8 }
  \}} \right).
\end{eqnarray}
 Uniformly on the support of $\pi$,
 $$\lim_{n \rightarrow \infty} b_n(f_0,f_i) = \frac{ 1 }{ 2 \pi }
\int_{-\pi}^\pi (1 - f_0/f_i)^2(x) dx = b(f_0,f_i),$$
$$ \lim_{n
\rightarrow \infty} h_n(f_0,f_i) = h(f_0,f_i).$$ Therefore, for
any $a>0$, if $n$ is large enough
\begin{eqnarray*} 
E_0^n[\phi_i] \leq   \max\left( \exp{\{-n\frac{ (h(f_0,f_i)^2-a)
}{ 16 (b(f_0,f_i)+ a) } \}},\exp{\{-n\frac{ h(f_0,f_i) + a }{ 8 }
  \}} \right),
\end{eqnarray*}
choosing $a < \varepsilon^2/2$, since $f_i \in A_\varepsilon^c$,
we obtain
\begin{eqnarray}   \label{test:null:1}
E_0^n[\phi_i] \leq   \max\left( \exp{\{-n\frac{ h(f_0,f_i)^2 }{ 32
(b(f_0,f_i)+\epsilon^2/2) } \}},\exp{\{-n\frac{ h(f_0,f_i)  }{ 16
} \}} \right).
\end{eqnarray}
Lemma \ref{ineq:b:h} implies that if $\varepsilon >0$ is small
enough,  there exists a constant $C_1$ such that
\begin{eqnarray*}   
E_0^n[\phi_i] \leq   \exp{(-n C_1 \epsilon
|\log{(\varepsilon)}|^{-1})}.
\end{eqnarray*}
 Moreover, if $f$ is in the support of $\Pi$ and satisfies $f\leq f_i$,
and $4|d_0-d| \leq \gamma$, using the same kind of calculations as
in the case of $E_0^n[\phi_i]$ and the fact that
 $$\id  - 2 s T_n^{1/2}(f)(T_n^{-1}(f_i) - T_n^{-1}(f_0))T_n^{1/2}
 \geq \id  + 2s B(f,f_0),$$
we obtain if $0<s<1/4$,
\begin{eqnarray*}
&E_f^n\left[ 1 - \phi_i \right] &\leq e^{ns\rho_i} \exp\{ - s
\tr{B(f,f_0)} + 4s^2 \tr{B(f,f_0)^2} \} \\
 &\leq& \exp\{ - ns h_n(f_0,f_i) + 4s^4 \tr{B(f,f_0)^2} +
 s\tr{A(f_i-f,f_0)} \}.
 \end{eqnarray*}
Using Lemma 2 if $(2\pi)^{-1}\int(f_i-f)f_0^{-1}(x)dx \leq
h(f_0,f_i)/4$, when $n$ is large enough (uniformly in $f$)
$\tr{A(f_i-f,f_0)} \leq n h_n(f_0,f_i)/2$ and
\begin{eqnarray} \label{test:alt:1}
E_f^n\left[ 1 - \phi_i \right] &\leq& \max \left( e^{-n\frac{
h(f_0,f_i)^2 }{ 32 b(f_i,f_0) }}, e^{-n \frac{ h(f_0,f_i) }{ 4 }
}\right).
\end{eqnarray}
Again Lemma \ref{ineq:b:h} implies that if $\varepsilon >0$ is
small enough, there exists a constant $C_1>0$ such that
\begin{eqnarray*}
E_f^n\left[ 1 - \phi_i \right] &\leq&  e^{-nC_1 \varepsilon
|\log{(\varepsilon)}|^{-1} }.
\end{eqnarray*}

\item If $4(d_i-d_0)> \gamma$. Set $\rho_i =
\tr{\id - T_n(f_0)T_n^{-1}(f_i)}/n + KL_n(f_0;f_i)$, the upper
bound of $E_0^n[\phi_i]$ is computed similarly to
(\ref{test:null:1}) so that
\begin{eqnarray*}
E_0^n[\phi_i] \leq   \max\left( \exp{\{-n\frac{ KL_n(f_0,f_i)^2 }{
8b_n(f_0,f_i) } \}},\exp{\{-n\frac{ KL_n(f_0,f_i) }{ 4 }
  \}} \right).
\end{eqnarray*}
Now, using the same calculations as in Dahlhaus (1989, p. 1754),
there exists a constant $C \geq 1$ such that
 $$KL_n(f_0,f_i) \geq \frac{ b_n(f_0,f_i) }{  C }$$
 so that, for large $n$ (independently of $f_i$),
 $$KL_n(f_0,f_i) \geq \frac{ b(f_0,f_i) }{ 2 C }$$
 We finally obtain that there exists a constant $c>0$ such that
\begin{eqnarray}   \label{test:null:2}
E_0^n[\phi_i] \leq    \exp{\{-n c b(f_0,f_i) \}}  .
\end{eqnarray}
Since $b(f_0,f_i)$ is uniformly bounded from below on the  set
$$\{
f \in \bar{\mathcal G}(t,m,M,L,\rho) ; d \geq  d_0 + 1/4(\rho
\wedge \rho_0 \wedge 1/2) - t/4 \}$$
(or $\bar{\mathcal
L}(t,m,M,L)$), if $\varepsilon $ is small enough
\begin{eqnarray*}
E_0^n[\phi_i] \leq    \exp{\{-n \varepsilon  \}} .
\end{eqnarray*}

Consider $f \leq f_i$, such that $4(d_i-d)\leq \rho \wedge 1/2
-t$. Similarly to before, let $h \in (0,1)$:
\begin{eqnarray*}   
&& E_f^n \left[ 1 - \phi_i\right] \leq  e^{(1-h)n\rho_i/2 - \frac{
1 }{ 2 } \log{\det[\id - (1-h)T_n^{1/2}(f)(T_n^{-1}(f_i) -
T_n^{-1}(f_0))T_n^{1/2}]}} \\
& \leq& e^{(1-h)n\rho_i/2 - \frac{ 1 }{ 2 } \log{\det[\id -
(1-h)B(f,f_0)]}} \\
&=& e^{(1-h)n\rho_i/2 - \log{\det[A(f,f_0)]}/2 - \frac{ 1 }{ 2 }
\log{\det[\id (1-h) - h T_n^{-1/2}(f)T_n(f_0)T_n^{-1/2}(f)]}},
 \end{eqnarray*}
then using the same kind of expansions as in (\ref{test:alt:1}),
we obtain
\begin{eqnarray*}   
&& E_f^n \left[ 1 - \phi_i\right] \leq   \det[A(f_i,f)]^{1/2}\\
&\times& \max\left\{ \exp\left( -
 \frac{ n KL_n(f_0,f_i)^2 }{ 32 \tr{B(f_0,f)^2}/n }\right), \exp\left( -
 \frac{ n KL_n(f_0,f_i) }{ 8 }\right) \right\}.
\end{eqnarray*}
Since $\log{\det[A(f,f_i)]} =-\log{\det[ \id +
T_n(f_i-f)T_n(f)^{-1}]}$, using a Taylor expansion of $\log{\det}$
around $\id$ , we obtain that for $n$ large enough
$$-\log{\det[A(f,f_i)]} \leq \frac{ 1 }{ 2\pi }\int_{-\pi}^\pi
(f_i-f)/f(x)dx + a$$ where $a$ can be chosen as small as
necessary. Also
$$ \frac{ KL_n(f_0,f_i)^2 }{ 32 \tr{B(f_0,f)^2}/n } \geq \frac{
 cb(f_0,f_i)^2-a }{ b(f_0,f) +a }.$$
 Since
 $$b(f_0,f) \leq 2b(f_0,f_i) + \frac1\pi\left( \int_{-\pi}^\pi f_0^2(f^{-1} -
 f_i^{-1})^2(x)dx\right),$$
 if
 $$ \frac{ 1 }{ 2\pi }\int_{-\pi}^\pi(f_i-f)/f(x)dx \leq
 cb(f_0,f_i) /4,$$
 $$\frac{ 1 }{ 2\pi} \int_{-\pi}^\pi f_0^2(f^{-1} -
 f_i^{-1})^2(x)dx \leq b(f_0,f_i),
 $$
 there exists a constant $c_1>0$ such that when $n$ is large
 enough
 \begin{eqnarray*}   
E_f^n \left[ 1 - \phi_i\right] &\leq & \exp\left( -nc_1b(f_0,f_i)
\right) \leq e^{-n\varepsilon}
\end{eqnarray*}
for $\varepsilon$ small enough.

\item If $4(d_0-d_i)>\gamma$. Set$\rho_i =
\log{\det[ T_n(f_i)T_n(f_0)^{-1}]}/n$, then if $0<h<1$
\begin{eqnarray*}
E_0^n \left[ \phi_i \right]
 &\leq&
 e^{-(1-h)n\rho_i/2 - \log{\det[A(f_0,f_i)]}/2- \log{\det[
 \id(1-h) + hT_n^{-1/2}(f_0)T_n(f_i)T_n^{-1/2}(f_0)]}} \\
  &\leq& e^{-nKL_n(f_i,f_0) + h^2 \tr{B(f_i,f_0)^2}} \leq e^{-n\varepsilon}
\end{eqnarray*}
where the last inequality can be obtained by following the same lines as
for (\ref{test:null:2}).

Moreover, for all  $f\leq f_i$, satisfying $4(d_i-d) \leq \rho
\wedge 1/2 - t$, if
 $$ \frac{ 1 }{ 2\pi }\int_{-\pi}^\pi(f_i-f)/f(x)dx \leq
 cb(f_i,f_0) /4, \quad  \frac{ 1 }{ 2\pi} \int_{-\pi}^\pi f_0^{-2}(f -
 f_i)^2(x)dx \leq b(f_i,f_0),
 $$
 where the constant $c$ is defined such that, for $n$ large
 enough, for all $f$ in $\bar{\mathcal G}(t,m,M,L,\rho)$ (resp. $\bar{\mathcal L}(t,m,M,L)$)
  such that $4(d_0-d)>\gamma$,
 $KL_n(f_i,f_0) \geq cb(f_i,f_0)$ (see the calculations presented in
 the case $4(d-d_0)>\gamma$),
\begin{eqnarray*}
  E_f^n \left[ 1 - \phi_i\right] &\leq & \exp\{ -2s n
  (KL_n(f_i,f_0) - \tr{ A(f_i-f,f)}/n ) + 4s^2 n b_n(f,f_0) \} \\
   &\leq& e^{-nc_1 b(f_i,f_0) }
\leq e^{-n\varepsilon}
\end{eqnarray*}
for $\varepsilon$ small enough.
\end{enumerate}
Then, in each case, we have, for large $n$
(independently of $f_i$),
 $$
 E_0^n \left[ \phi_i \right] \leq e^{-n \varepsilon |\log \varepsilon|^{-1}}
$$
for all $\varepsilon < \varepsilon_0$,
$$E_f^n \left[ 1 - \phi_i
\right] \leq e^{-n \varepsilon |\log \varepsilon|^{-1}} .$$
Let
$\phi^{(n)} = \max_i \phi_i$; then, using Markov inequality,
\begin{eqnarray*}
\tilde{p}_2 &\leq& E_0^n \left[ \phi^n \right] +
 2e^{2n\delta}\int_{A_\varepsilon
\cap {\cal F}_n} E_f \left[ 1 - \phi^n \right] d\pi(f) \\
 &\leq& e^{N_n}e^{-n\varepsilon|\log \varepsilon|^{-1}}  + 2e^{2n\delta}e^{-n\varepsilon|\log \varepsilon|^{-1} } \\
 &\leq& e^{-n\varepsilon |\log \varepsilon|^{-1}/2 },
 \end{eqnarray*}
We finally obtain that for some $\delta >0$, if $n$ is large enough
 $$P_0^n\left[ P^\pi[A_\varepsilon^c|\Xn ] > e^{-n\delta}  \right] \leq
 \frac{ C_0 }{ n^3} $$
for some positive constant $C_0$, so that $\pi[A_\varepsilon^c|\Xn]
\rightarrow 0 $ $P_0^\infty$ a.s.

\section{\texorpdfstring{Lemma \ref{whittles}}{Lemma 3}} \label{app:lemmas}

\begin{lemma}\label{whittles}
Let $f_j$, $j \in \{1,2\}$ be such that $f_j(\lambda) =
|\lambda|^{-2 d_j} \tilde{f}_j(\lambda)$,
 where $d_j < 1/2$ and $\tilde{f}_j \in \mathcal{S}(L,\beta)$, for some constant
 $L>0$ and consider $b$ a bounded function on $[-\pi, \pi]$.
Assume that $h_n(f_1,f_2) < \epsilon$ where $\epsilon> 0$.  Then
$\forall \delta > 0$, there exists $\epsilon_0> 0$ such that if
$\epsilon < \epsilon_0$, there exists $C>0$ such that
\beq\label{whit0}
\frac{ 1 }{ n  } \tr{ T_n(f_1)^{-1}T_n(f_1b)T_n(f_2)^{-1}T_n(f_1b) }
\leq C (\log{n})^3 [ |b|_2^2 + |b|_\infty^2
  n^{\delta-1} + n^{-1/2+\delta}],
 \enq
\begin{eqnarray}\label{whit01}
\frac{ 1 }{ n  } \tr{
T_n(f_1^{-1})T_n(f_1-f_2)T_n(f_2^{-1})T_n(f_1 - f_2) } &\leq&
C h_n(f_1,f_2).
 \end{eqnarray}
Let $g_j = (1 - \cos{\lambda})^{d_j}$ and $f_j =
g_j^{-1}\tilde{f}_j$, where $\tilde{f}_1 \in \mathcal{S}(L,\beta)
\cap \mathcal{L}$ and $\tilde{f}_2 \in \mathcal{S}(L,\beta)$,
written in the form
$\log\tilde{f}_2(\lambda) = \sum_{l=0}^K \theta_l \cos{(l \lambda)}$; then
\begin{eqnarray}\label{whit02}
 & & \left|\frac{ 1 }{ n }\tr{ T_n(g_1(f_1-f_2))T_n(g_2(f_1-f_2))} -
 \tr{T_n(g_1g_2(f_1-f_2)^2)} \right| \nonumber \\
 & \leq & Cn^{-1+\delta} + n^{-1}\log{n} \sum_{l=0}^{K_n} l|\theta_l|
\left(
  \int_{[-\pi,\pi]}g_1g_2(f_1-f_2)^2(\lambda)d\lambda\right)^{1/2},
\end{eqnarray}
for any $\delta >0$.
\end{lemma}

\begin{proof}
Throughout the proof $C$ denotes a generic constant. We first
prove (\ref{whit0}). To do so, we obtain an upper bound on another
quantity, namely
\begin{eqnarray}
 \gamma(b) &=&\frac{ 1 }{ n }\tr{
 T_n(f_1^{-1})T_n(f_1b)T_n(f_2^{-1})T_n(f_1b)}.
\end{eqnarray}
 First note that $b$ can be replaced by
$|b|$ so that we can assume that it is positive. Let
$\Delta_n(\lambda) = \sum_{j=1}^n\exp( -i\lambda j)$ and $L_n$ be
the $2\pi$-periodic function defined by $L_n(\lambda) = n$ if
$|\lambda| \leq 1/n$ and $L_n(\lambda) = |\lambda|^{-1}$ if $1/n
\leq |\lambda| \leq \pi$. Then $|\Delta_n(\lambda)| \leq
CL_n(\lambda)$  and we can express traces of products of Toeplitz
matrices in the following way. Let the symbol $d\bla$ denote the quantity
$d\la_1 d\la_2 d\la_3 d\la_4$;
\begin{eqnarray*}
 \gamma(b) &=&   \frac{ C }{ n }\int_{[-\pi,\pi]^4} \frac{
  f_1(\lambda_1)b(\lambda_1)f_1(\lambda_3)b(\lambda_3)}{ f_0(\lambda_2)f_0(\lambda_4) } \times \\
   & & \quad
\Delta_n(\lambda_1-\lambda_2)\Delta_n(\lambda_2-\lambda_3)\Delta_n(\lambda_3-\lambda_4)\Delta_n(\lambda_4-\lambda_1)d\bla\\
   &=&
\frac{ C }{ n }\int_{[-\pi,\pi]^4}  \frac{
  f_1(\lambda_1)f_1(\lambda_3)}{ f_0(\lambda_2)f_0(\lambda_4) }  \left( b(\lambda_1)^2 +
  b(\lambda_1)b(\lambda_3) - b(\lambda_1)^2 \right) \\
  & & \quad \times
  \Delta_n(\lambda_1-\lambda_2)\Delta_n(\lambda_2-\lambda_3)\Delta_n(\lambda_3-\lambda_4)\Delta_n(\lambda_4-\lambda_1)d\bla\\
  &=& \frac{ C }{ n }\tr{
 T_n(f_1b^2)T_n(f_1^{-1})T_n(f_1)T_n(f_2^{-1})} \\
  & & +
\frac{ C }{ n }\int_{[-\pi,\pi]^4} \frac{
  f_1(\lambda_1)f_1(\lambda_3)b(\lambda_1) }{ f_1(\lambda_2)f_2(\lambda_4) }\left( b(\lambda_3) -
  b(\lambda_1) \right) \\
  & & \quad \times
\Delta_n(\lambda_1-\lambda_2)\Delta_n(\lambda_2-\lambda_3)\Delta_n(\lambda_3-\lambda_4)\Delta_n(\lambda_4-\lambda_1)d\bla
 \end{eqnarray*}
On the set $b(\lambda_1)> b(\lambda_3)$, $0 < b(\lambda_1) -
b(\lambda_3)< b(\lambda_1)$ and on the set $b(\lambda_3)>
b(\lambda_1)$, $0 < b(\lambda_3) - b(\lambda_1)< b(\lambda_3)$,
therefore the second term of the r.h.s. of the above inequality is
bounded by (in absolute value)
 \begin{eqnarray*}
 \gamma(b)  &\leq & \frac{ C }{ n }
 \int_{[-\pi,\pi]^4} \frac{
  f_1(\lambda_1)f_1(\lambda_3)b(\lambda_1)^2 }{ f_1(\lambda_2)f_2(\lambda_4) } L_n(\lambda_1-\lambda_2)L_n(\lambda_2-\lambda_3) \\
  & & \quad \quad \times
 L_n(\lambda_3-\lambda_4)L_n(\lambda_4-\lambda_1)d\bla \\
 &\leq& \frac{ C}{ n }\int_{[-\pi,\pi]^4} b(\lambda_1)^2\frac{
 |\lambda_1|^{-2d_1}|\lambda_3|^{-2d_1} }{
 |\lambda_2|^{-2d_1}|\lambda_4|^{-2d_2} }
L_n(\lambda_1-\lambda_2)L_n(\lambda_2-\lambda_3) \\
  & & \quad \quad \times
 L_n(\lambda_3-\lambda_4)L_n(\lambda_4-\lambda_1)d\bla
\end{eqnarray*}
Note that
\begin{eqnarray} \label{ineq:int:Ln}
\int_{[-\pi,\pi]}L_n(\lambda_1-\lambda_2)L_n(\lambda_2-\lambda_3)d\lambda_2
\leq C\log{n}L_n(\lambda_1-\lambda_3), \end{eqnarray}
 therefore
  \begin{eqnarray*}
 \gamma(b)  & \leq &
 \frac{ C(\log{n})^3}{ n }\int_{[-\pi,\pi]}b(\lambda)^2d\lambda \\
  & &    +   C \int_{[-\pi,\pi]^4}b(\lambda_1)^2|\lambda_1|^{-2(d_1-d_2)}\left(  \frac{
|\lambda_3|^{-2d_1} }{ |\lambda_2|^{-2d_1} } - 1\right)\left(
\frac{
|\lambda_1|^{-2d_2} }{ |\lambda_4|^{-2d_2} }  -1 \right) \\
  & & \quad \times
  L_n(\lambda_1-\lambda_2)L_n(\lambda_2-\lambda_3)L_n(\lambda_3-\lambda_4)L_n(\lambda_4-\lambda_1)d\bla \\
    & &
    + 2C \int_{[-\pi,\pi]^4}b(\lambda_1)^2|\lambda_1|^{-2(d_1-d_2)}\left(  \frac{
|\lambda_3|^{-2d_1} }{ |\lambda_2|^{-2d_1} } +  \frac{
|\lambda_1|^{-2d_2} }{ |\lambda_4|^{-2d_2} }
 -2  \right) \\
  & & \quad \times
  L_n(\lambda_1-\lambda_2)L_n(\lambda_2-\lambda_3)L_n(\lambda_3-\lambda_4)L_n(\lambda_4-\lambda_1)d\bla
\end{eqnarray*}

Since \begin{eqnarray}\label{ineq:Dahl}
 \left| \frac{
|\lambda_1|^{-2d_j} }{ |\lambda_2|^{-2d_j} } - 1 \right| \leq
  C \frac{ |\lambda_1 - \lambda_2 |^{1-\delta} }{ |\lambda_1|^{1-\delta}
  },\quad \mbox{for } j=\{1,2\},
\end{eqnarray}
  using Dahlhaus' (1989) calculations as in his proof of Lemma 5.2, we obtain
  that, if $d_1 - d_2 < \delta/4$,
\begin{eqnarray*}
\lefteqn{
\int\limits_{[-\pi,\pi]^4}b(\lambda_1)^2|\lambda_1|^{-2(d_1-d_2)}\left(
\frac{ |\lambda_3|^{-2d_1} }{ |\lambda_2|^{-2d_1} } -
1\right)\left( \frac{ |\lambda_1|^{-2d_2} }
 { |\lambda_4|^{-2d_2} }-1  \right) } \\
  &\times &
   L_n(\lambda_1-\lambda_2)L_n(\lambda_2-\lambda_3)L_n(\lambda_3-\lambda_4)
   L_n(\lambda_4-\lambda_1)d\bla \\
   &\leq& |b|_\infty^2\int\limits_{[-\pi,\pi]^4}\frac{
L_n(\lambda_1\!-\!\lambda_2)L_n(\lambda_2\!-\!\lambda_3)^{\delta}L_n(\lambda_3\!-\!\lambda_4)L_n(\lambda_4\!-\!\lambda_1)^{\delta}}{|\lambda_1|^{ 1 -\delta/2} |\lambda_4|^{1 -\delta}}
d\bla\\
  &\leq& C n^{2\delta}|b|_\infty^2(\log{n})^2,
\end{eqnarray*}
as long as $|d_1 - d_2| < \delta/2$.
 By considering
$h_n(f,f_0)<\epsilon$ with $\epsilon>0$ small enough, we can
impose that $|d_1-d_2|< \delta/2$, and we finally obtain that
 \begin{eqnarray} \label{res:f:1}
  \gamma(b) &\leq&  C |b|_2^2 (\log{n})^3 + C|b|_\infty^2
  n^{2\delta-1}(\log{n})^2.
\end{eqnarray}
We now prove  that, for large $n$ and $\forall \delta>0$,
 $$ \frac{ 1 }{ n  } \tr{
T_n(f_1)^{-1}T_n(f_1b)T_n(f_2)^{-1}T_n(f_1b) } \leq C \frac{ 1 }{ n
} \tr{ T_n(f_1^{-1})T_n(f_1b)T_n(f_2^{-1})T_n(f_1b) } +
Cn^{-1+\delta}.$$
Since
 $f_i(\lambda) \leq C |\lambda|^{-2d_i} \propto g_i(\lambda)$, $i=1,2$. This implies that
$T_n^{-1}(f_i) \succeq C^{-1} T_n^{-1}(g_i)$ so that we can replace
$T_n(f_i)^{-1}$ by $T_n^{-1}(g_i)$ in the above term. Then
\begin{eqnarray*}
 \delta_n &=& \tr{T_n(f_1b)T_n^{-1}(g_1)T_n(f_1b)T_n^{-1}(g_2)} \\
  &=&
 \tr{T_n(f_1b)T_n(g_1^{-1}/(4\pi^2))T_n(f_1b_n)
 T_n(g_2^{-1}/(4\pi^2)))} \\
 &+&  \tr{T_n(f_1b)T_n^{-1}(g_1)T_n(f_1b)
 T_n^{-1/2}(g_2)R_2T_n^{-1/2}(g_2)} \\
  &+&
  \tr{T_n(f_1b)T_n(g_1)^{-1/2}R_1T_n(g_1)^{-1/2}T_n(f_1b) T_n(g_2^{-1}/(4\pi^2))},
\end{eqnarray*}
where $R_i = T_n(g_i)^{1/2}T_n(g_i^{-1}/(4\pi^2))T_n(g_i)^{1/2}-
\id$, $i=1,2$. Using (\ref{res:f:1}) we obtain that
\begin{eqnarray*}
\tr{ T_n(f_1b)T_n(g_1^{-1})T_n(f_1b) T_n(g_2^{-1}) } &\leq &
C(\log{n})^3  n|b_n|_2^2 + O(|b|_\infty n^{\delta}) = n\gamma.
\end{eqnarray*}
Moreover
\begin{eqnarray*}
\lefteqn{\left|\tr{T_n(f_1b)T_n^{-1}(g_1)T_n(f_1b)
T_n^{-1/2}(g_2)R_2T_n^{-1/2}(g_2)} \right|   } \\
 & \leq& |R_2||
T_n^{-1/2}(g_2)T_n(f_1b)T_n^{-1}(g_1)T_n(f_1b) T_n^{-1/2}(g_2)| \\
 &\leq&\delta_n^{1/2} |R_2|  |\!|T_n^{-1/2}(g_2)T_n(f_1b)^{1/2}|\!| |\!|T_n(f_1b)^{1/2}
  T_n^{-1/2}(g_1)|\!|
\end{eqnarray*}
Lemmas 5.2 and 5.3 in Dahlhaus (1989) lead to: $\forall \delta >0$
\begin{eqnarray*}
\left|\tr{T_n(f_1b)T_n^{-1}(g_1)T_n(f_1b)
T_n^{-1/2}(g_2)R_2T_n^{-1/2}(g_2)} \right|  \leq Cn^{\delta +
2(d_0-d)}|b_n|_\infty \delta_n^{1/2} \leq C \delta_n^{1/2}
\end{eqnarray*}
Similarly,
\begin{eqnarray*}
&\,& \left|\tr{T_n(f_1b)T_n(g_1)^{-1/2}R_1T_n(g_1)^{-1/2}T_n(f_1b)
T_n(g_1^{-1}/(4\pi^2))} \right|\\
&\leq& |R_1|\delta_n^{1/2} [|R_1|+ 1]
|\!|T_n(g_1)^{-1/2}T_n(f_1b)^{1/2}|\!|
 |\!|
T_n(f_1b)^{1/2}T_n(g_1)^{-1/2}|\!|
\end{eqnarray*}
Since
$|\!|T_n(g_1)^{1/2}T_n(g_1^{-1}/(4\pi^2))^{1/2}|\!| \leq 1 +
|\!|T_n(g_1)^{1/2}T_n(g_1^{-1}/(4\pi^2))^{1/2}- \id|\!| \leq
n^\delta$ for all $\delta >0$ and using Lemma 5.3 of Dahlhaus (1989)
\begin{eqnarray*}
\left|\tr{T_n(f_1b)T_n(g_1)^{-1/2}R_1T_n(g_1)^{-1/2}T_n(f_1b)
T_n(g_2^{-1}/(4\pi^2))} \right|&\leq& C n^{2\delta}\delta_n^{1/2}.
\end{eqnarray*}
Finally we obtain for all $\delta >0$, when $n$ is large enough
\begin{eqnarray*}
\delta_n/n \leq Cn^{-1/2+\delta}\sqrt{\delta_n/n} + \gamma/n \leq 2
\gamma/n + 0(n^{-1+\delta}),
\end{eqnarray*}
 and (\ref{whit0}) is proved. We now prove (\ref{whit01}).
since $f_j \geq m|\lambda|^{-2d_j} = g_j$ where $m = e^{-L}$,
 $T_n^{-1}(f_j) \prec T_n^{-1}(g_j)$, i.e. $T_n^{-1}(g_j)-T_n^{-1}(f_j)$ is positive semidefinite, and
\begin{eqnarray}
&& h_n(f_1,f_2) =  \frac{ 1 }{ 2n }\tr{T_n(f_1-f_2)
T_n^{-1}(f_2)T_n(f_1-f_2)T_n^{-1}(f_1) } \nonumber \\
 &\geq& \frac{ 1 }{ 2n }\tr{T_n(f_1-f_2) T_n^{-1}(f_2)T_n(f_1-f_2)T_n^{-1}(g_1)} \nonumber \\
 &\geq&
\frac{ 1 }{ 2n }\tr{T_n(f_1-f_2)
  T_n^{-1}(f_2)T_n(f_1-f_2)T_n^{-1/2}(g_1)R_1T_n^{-1/2}(g_1)} \nonumber \\
\end{eqnarray}
\begin{eqnarray} \label{ineq1:whit01}
&+&
 \frac{ 1 }{ 2n }\tr{T_n(f_1-f_2) T_n^{-1}(g_2)T_n(f_1-f_2)T_n\left(\frac{ g_1^{-1}}{ 4\pi^2 } \right)
} \nonumber \\
  &= & \frac{ 1 }{ 2n(16\pi^4) }\tr{T_n(f_1-f_2) T_n(g_2^{-1})T_n(f_1-f_2)T_n\left(g_1^{-1}\right )
} \nonumber \\
 &+& \frac{ 1 }{ 2n }\tr{T_n(f_1-f_2)
  T_n^{-1}(f_2)T_n(f_1-f_2)T_n^{-1/2}(g_1)R_1T_n^{-1/2}(g_1)}
   \\
  &+& \frac{ 1 }{ 2n(4\pi^2) }\tr{T_n(f_1-f_2) T_n^{-1/2}(g_2) R_2T_n^{-1/2}(g_2)T_n(f_1-f_2)T_n\left( g_1^{-1}\right)} \nonumber
 \end{eqnarray}
where $R_j = \id -
T_n^{1/2}(g_j)T_n(g_j^{-1}/(4\pi^2))T_n^{1/2}(g_j)$. We first bound
the first term of the r.h.s. of (\ref{ineq1:whit01}).
Let $\delta > 0$ and $\epsilon< \epsilon_0$ such that $|d-d_0| \leq \delta$
(Corollary \ref{Cor:cons:d} implies that there exists such a value
$\epsilon_0$). Then using Lemmas 5.2 and 5.3 of
\cite{dhl:89}
\begin{eqnarray*}
\lefteqn{ \left|\tr{T_n(f_1-f_2)
  T_n^{-1}(f_2)T_n(f_1-f_2)T_n^{-1/2}(g_1)R_1T_n^{-1/2}(g_1)}
  \right|} \\
  &\leq& 2
 |R_1||T_n^{-1/2}(g_1)T_n(f_1\!\!-\!\!f_2)T_n^{-1/2}(f_2)|
 |\!|T_n(|f_1\!\!-\!\!f_2|)^{1/2}T_n^{-1/2}(f_2)|\!| \\
 &\times&
 |\!|T_n(|f_1\!\!-\!\!f_2|)^{1/2}T_n^{-1/2}(g_1)|\!| \\
  &\leq& Cn^{3\delta}
  |T_n^{-1/2}(g_1)T_n(f_1\!\!-\!\!f_2)T_n^{-1/2}(f_2)|.
\end{eqnarray*}
Since $g_1 \leq C f_1$,
 \begin{eqnarray*}
 |T_n^{-1/2}(g_1)T_n(f_1\!\!-\!\!f_2)T_n^{-1/2}(f_2)|^2 &=&
 \tr{T_n^{-1}(g_1)T_n(f_1\!\!-\!\!f_2)T_n^{-1}(f_2)T_n(f_1\!\!-\!\!f_2)} \\
  &\leq& C \tr{T_n^{-1}(f_1)T_n(f_1\!\!-\!\!f_2)T_n^{-1}(f_2)T_n(f_1\!\!-\!\!f_2)} \\
   &=& C n h_n(f_1,f_2),
 \end{eqnarray*}
 and
\begin{eqnarray*}
\frac{ 1 }{ n }\left|\tr{T_n(f_1\!\!-\!\!f_2)
  T_n^{-1}(f_2)T_n(f_1\!\!-\!\!f_2)T_n^{-1/2}(g_1)R_1T_n^{-1/2}(g_1)}
\right| &\leq&  C n^{2\delta-1/2}h_n(f_1,f_2).
\end{eqnarray*}
We now bound the second term of the r.h.s. of (\ref{ineq1:whit01}).
\begin{eqnarray*}
   &= & \left| \frac{ 1 }{ n }\tr{T_n(f_1\!\!-\!\!f_2) T_n^{-1/2}(g_2)R_2T_n^{-1/2}(g_2)T_n(f_1\!\!-\!\!f_2)T_n(g_1^{-1})}  \right| \\
 &\leq&
 \frac{ 1 }{ n } |R_2| |T_n^{-1/2}(g_2)T_n(f_1\!\! -
 \!\!f_2)T_n(g_1)^{-1/2}| |
 T_n(g_1)^{1/2}T_n(g_1^{-1})T_n(|f_1\!\!-\!\!f_2|) T_n^{-1/2}(f_2)| \\
  & \leq & \frac{ Cn^{\delta} \sqrt{nh_n(f_2,f_1)} }{ n }
  |\!| T_n(g_1)^{1/2}T_n(g_1^{-1})T_n(|f_1\!\!-\!\!f_2|) T_n^{-1/2}(f_2)|\!| \\
  &\leq& \frac{ Cn^{\delta+1/2} \sqrt{h_n(f_2,f_1)} }{ n }
  |\!| T_n(g_1)^{1/2}T_n(g_1^{-1})^{1/2}|\!|^2  \\
    & & \times |\!|T_n(g_1)^{-1/2}T_n(|f_1\!\!-\!\!f_2|)^{1/2}|\!| |\!| T_n(|f_1\!\!-\!\!f_2|)^{1/2} T_n^{-1/2}(f_2)|\!|\\
  &\leq & C n^{3\delta-1/2}h_n(f_1,f_2),
\end{eqnarray*}
since  $|\!|T_n(f_1)^{1/2}T_n(f_1^{-1})T_n(f_1)^{1/2}|\!| \leq
|\!| \id |\!| + |T_n(f_1)^{1/2}T_n(f_1^{-1})T_n(f_1)^{1/2}-\id|
\leq C n^\delta$.\\

\noindent
Therefore,
$$
  \frac{ C }{ n }
  \tr{T_n(f_1\!\!-\!\!f_2)
T_n(g_2^{-1})T_n(f_1\!\!-\!\!f_2)T_n(g_1^{-1})} \leq
 C\,h_n(f_1,f_2) ( 1 +n^{-1/2 + 3 \delta}),
$$
and, using the fact that $C\, g_j > f_j$, for $j=1,2$ this
proves (\ref{whit01}).
The proof of (\ref{whit02}) is similar:
\begin{eqnarray*}
 A &=& \tr{ T_n(g_1(f_1-f_2))T_n(g_2(f_1-f_2))} -
 \tr{
T_n(g_1g_2(f_1-f_2)^2)} \\
 &=&
C\int\limits_{[-\pi,\pi]^2}\!\!\!g_1(f_1\!-\!f_2)(\lambda_1)[g_2(f_1\!-\!f_2)(\lambda_2) - g_1(f_1\!-\!f_2)(\lambda_1)
 ] \Delta_n(\lambda_1\!-\!\lambda_2)...\Delta_n(\lambda_4\!-\!\lambda_1)d\bla \\
 &=& C\int\limits_{[-\pi,\pi]^2}\!\!g_1(f_1-f_2)(\lambda_1)(f_1-f_2)(\lambda_2)[g_2(\lambda_2) -
 g_2(\lambda_1)]
 \Delta_n(\lambda_1\!-\!\lambda_2)\Delta_n(\lambda_2\!-\!\lambda_1)d\bla\\
  &-&  C\int\limits_{[-\pi,\pi]^2}g_1g_2(f_1-f_2)(\lambda_1)[f_1(\lambda_2)-f_1(\lambda_1)]
 \Delta_n(\lambda_1-\lambda_2)\Delta_n(\lambda_2-\lambda_1)d\bla\\
  &+&    C\int\limits_{[-\pi,\pi]^2}g_1g_2(f_1-f_2)(\lambda_1)[f_2(\lambda_2)-f_2(\lambda_1)]
 \Delta_n(\lambda_1-\lambda_2)\Delta_n(\lambda_2-\lambda_1)d\bla.
\end{eqnarray*}
The first 2 terms of the right hand side are of order
$O(n^{2\delta}\log{n})$. We now study the last term, here the
problem is due to the fact that $\tilde{f}_2$ does not necessarily
belong to $\mathcal{L}$. We have:
\begin{eqnarray*}
&\int&_{[-\pi,\pi]^2}g_1g_2(f_1-f_2)(\lambda_1)[f_2(\lambda_2)-f_2(\lambda_1)]
 \Delta_n(\lambda_1-\lambda_2)\Delta_n(\lambda_2-\lambda_1)d\bla
\\
&=&
\int\limits_{[-\pi,\pi]^2}g_1g_2(f_1-f_2)(\lambda_1)\tilde{f}_2(\lambda_2)[g_2^{-1}(\lambda_2)-g_2^{-1}(\lambda_1)]
 \Delta_n(\lambda_1-\lambda_2)\Delta_n(\lambda_2-\lambda_1)d\bla
 \\
&+&
\int\limits_{[-\pi,\pi]^2}g_1(f_1-f_2)(\lambda_1)[\tilde{f}_2(\lambda_2)-\tilde{f}_2(\lambda_1)]
 \Delta_n(\lambda_1-\lambda_2)\Delta_n(\lambda_2-\lambda_1)d\bla.
\end{eqnarray*}
The first term of the above inequality is of order
$O(n^{2\delta}\log{n})$ because $g_2$ belongs to $\mathcal{L}.$
Since
 $$\tilde{f}(\lambda) = \exp\left (\sum_{l=0}^{K_n}\theta_l
 \cos{(l\lambda})\right ),$$ one gets
\begin{eqnarray*}
 I&=& \int\limits_{[-\pi,\pi]^2}g_1(f_1-f_2)(\lambda_1)[\tilde{f}_2(\lambda_2)-
 \tilde{f}_2(\lambda_1)]
 \Delta_n(\lambda_1-\lambda_2)\Delta_n(\lambda_2-\lambda_1)d\bla
 \\
 &\leq&
 C \int\limits_{[-\pi,\pi]^2}g_1|f_1-f_2|(\lambda_1)\left|\sum_{j=0}^{K_n}
 \theta_l
 (\cos{(j \lambda_2)}-\cos{(j \lambda_1)})\right|
L_n(\lambda_1-\lambda_2)L_n(\lambda_2-\lambda_1)d\bla \\
 &\leq&
C \log{n} \left( \sum_{l=0}^{K_n}| \theta_l|l\right)
 \int\limits_{[-\pi,\pi]}g_1|f_1-f_2|(\lambda)d\bla \\
 &\leq& C \log{n} \sum_{l=0}^{K_n}| \theta_l|l \left(
 \int\limits_{[-\pi,\pi]}g_1g_2(f_1-f_2)^2(\lambda)d\bla\right)^{1/2},
\end{eqnarray*}
where the latter inequality holds
because $\int g_1/g_2(\lambda)d\bla$
can be proved to be bounded by an application of
an application of H\"older inequality. \end{proof}

\section{Relations between $b(f_0,f)$ and $h(f_0,f)$}

\begin{lemma} \label{ineq:b:h}
Let $m,M,L>0$ and $\rho \leq 1$. There exists  $\tau>0$  and $C>0$
such that for any $f,f_0 \in \bar{\mathcal G}(t,m,M,L,\rho) \cup
\bar{\mathcal L}(t,m,M,L)$, if $h(f,f_0)<\tau$,
$$b(f,f_0) \leq h(f,f_0) |\log{h(f,f_0)}|.$$
\end{lemma}
We need to bound $b(f,f_0)$ in terms of $h(f_0,f)$ when $|d-d_0|$ is
small.
Assume that $f_0 = |x|^{-2d_0}\tilde{f}_0$ and  $f =
|x|^{-2d}\tilde{f}$ with $d\geq d_0$ (otherwise the bound is
straightforward) we have
$$ b(f,f_0) = \frac{ 1 }{ 2\pi }\int\limits_{-\pi}^\pi (f/f_0 - 1)^2dx = \frac{ 1 }{ 2\pi }\int\limits_{-\pi}^\pi \frac{ (f-f_0)^2 }{ f_0^2 }dx,$$
$$ h(f,f_0) = \frac{ 1 }{ 2\pi }\int\limits_{-\pi}^\pi (f/f_0 -
1)^2\frac{ f_0}{f}dx= \frac{ 1 }{ 2\pi }\int\limits_{-\pi}^\pi \frac{
(f-f_0)^2 }{ f_0 f }dx.$$
Let $A>0$ be large enough; using the fact
that $m \leq \tilde f, \tilde f_0 \leq M$ we obtain, if $A >2$, that
\begin{eqnarray*}
b(f,f_0) &\leq & A h(f_0,f) + \int\limits_{f/f_0 > A}\frac{ 1 }{ 2\pi
}\int\limits_{-\pi}^\pi (f/f_0 - 1)^2dx \\
 &\leq& A h(f_0,f) + \int\limits_{f/f_0 > A}\frac{ 2M^2 }{ 2m^2\pi
} |x|^{-4(d-d_0)}dx.
\end{eqnarray*}
Let $A > M/m$ then if $f/f_0 > A $, $|x|^{-2(d-d_0)} > Am/M ,$
 so that
\begin{eqnarray*}
b(f,f_0) &\leq & A h(f_0,f) + C\int\limits_{|x|^{-2(d-d_0)} > K A}
|x|^{-4(d-d_0)}dx
\end{eqnarray*}
Now assume that $h(f_0,f) \leq \tau$ where $\tau >0$ is fixed and
small. Consider $t>0$ small enough so that
$$h(f_0,f) \gtrsim \int\limits_{x^{-2(d-d_0)} \geq t^{-1}}x^{-2(d-d_0)}dx$$
where $\gtrsim$ means that the inequality is up to a multiplicative
constant whose value does not depend on $f$ and $f_0$ (but it does depend on $M$
and $m$). It implies that
\begin{eqnarray*}
\frac{ 1 }{ 1 -2(d-d_0)}t^{- 1 +1/2(d-d_0)} \leq C h(f,f_0),
\end{eqnarray*}
so that, if $t^{-1} = \log{1/h(f,f_0)}$,
$$\log(\log{1/h(f,f_0)})\frac{ 1 }{  2(d-d_0) } \geq \log{1/(\rho h(f,f_0))}, \quad \rho >0 \mbox{ fixed }.$$
Hence, if $h(f,f_0)$ is small enough,
$$
2(d -d_0) \leq \frac{ 2\log(\log{1/h(f,f_0)}) }
{\log{1/h(f_0,f)}}.
$$
Now using the fact that there exists $C>0$ such that
 $$
b(f,f_0) \leq A h(f_0,f) +  C A^{2-1/(2(d-d_0))}, \quad h(f_0,f)
\geq C' \frac{  A^{1-2(d-d_0)} }{ 1 - 2(d-d_0)}$$ and considering
 $A = \log{1/h(f,f_0)}$ we finally obtain
 $$b(f,f_0) \leq  \log{1/h(f,f_0)} h(f_0,f) + C' h(f,f_0)
\log{(1/h(f,f_0))}.$$
Hence, there exists $\tau>0$ (depending only on $m,M$ and $C>0$)
such that if $h(f,f_0)<\tau$,
$b(f,f_0) \leq h(f,f_0) |\log{h(f,f_0)}|$.


\begin{thebibliography}{27}

\bibitem{aitc:80}
Aitchison J. and Shen S.M. (1980) Logistic-normal distributions: some
  properties and uses, {\em Biometrika\/}, 67, 2, 261--272.

\bibitem{bardet:03}
Bardet J.M., Lang G., Oppenheim G., Philippe A., Stoev S. and Taqqu M.S. (2003)
  Semi-parametric estimation of the long-range dependence parameter: a survey,
  in: {\em Theory and applications of long-range dependence\/}, Birkh\"auser,
  Boston, MA, 557--577.

\bibitem{ber:93}
Beran J. (1993) Fitting long-memory models by generalized linear regression,
  {\em Biometrika\/}, 80, 4, 817--822.

\bibitem{ber:94}
Beran J. (1994) {\em Statistics for long-memory processes\/}, volume~61 of {\em
  Monographs on Statistics and Applied Probability\/}, Chapman and Hall, New
  York.

\bibitem{blom:73}
Bloomfield P. (1973) An exponential model for the spectrum of a scalar time
  series, {\em Biometrika\/}, 60, 217--226.

\bibitem{Hur:2006}
Chen W.W., Hurvich C.M. and Lu Y. (2006) On the {C}orrelation {M}atrix of the
  {D}iscrete {F}ourier {T}ransform and the {F}ast {S}olution of {L}arge
  {T}oeplitz {S}ystem for {L}ong {M}emory {T}ime {S}eries, {\em J. Amer.
  Statist. Assoc.\/}, 101, 474, 812--821.

\bibitem{chgh:04}
Choudhuri N., Ghosal S. and Roy A. (2004) Bayesian estimation of the spectral
  density of a time series, {\em J. Amer. Statist. Assoc.\/}, 99, 468,
  1050--1059.

\bibitem{dhl:89}
Dahlhaus R. (1989) Efficient parameter estimation for self-similar processes,
  {\em Ann. Statist.\/}, 17, 4, 1749--1766.

\bibitem{douc:05}
Douc R., Guillin A., Marin J. and Robert C. (2006) Convergence of adaptive
  mixtures of importance sampling schemes, {\em Ann. Statist. (to appear)\/}.

\bibitem{dou:03}
Doukhan P., Oppenheim G. and Taqqu M.S. (Eds.) (2003) {\em Theory and
  applications of long-range dependence\/}, Birkh\"auser Boston Inc., Boston,
  MA.

\bibitem{ft:86}
Fox R. and Taqqu M.S. (1986) Large-sample properties of parameter estimates for
  strongly dependent stationary {G}aussian time series, {\em Ann. Statist.\/},
  14, 2, 517--532.

\bibitem{gewe:83}
Geweke J. and Porter-Hudak S. (1983) The estimation and application of long
  memory time series models, {\em J. Time Ser. Anal.\/}, 4, 4, 221--238.

\bibitem{ggvdv:01}
Ghosal S., Ghosh J.K. and van~der Vaart A.W. (2000) Convergence rates of
  posterior distributions, {\em Ann. Statist.\/}, 28, 2, 500--531.

\bibitem{gvdv:06}
Ghosal S. and Van~der Vaart A. (2007) Convergence rates of posterior
  distributions for non i.i.d. observations, {\em Ann. Statist.\/}, 35,
  192--225.

\bibitem{ghrm:00}
Ghosh J. and Ramamoorthi R. (2003) {\em Bayesian nonparametrics\/}, Springer
  Series in Statistics, Springer-Verlag, New York.

\bibitem{gt:99}
Giraitis L. and Taqqu M.S. (1999) Whittle estimator for finite-variance
  non-{G}aussian time series with long memory, {\em Ann. Statist.\/}, 27, 1,
  178--203.

\bibitem{hms:02}
Hurvich C.M., Moulines E. and Soulier P. (2002) The {FEXP} estimator for
  potentially non-stationary linear time series, {\em Stochastic Process.
  Appl.\/}, 97, 2, 307--340.

\bibitem{jud:03}
Lieberman O., Rousseau J. and Zucker D.M. (2003) Valid asymptotic expansions
  for the maximum likelihood estimator of the parameter of a stationary,
  {G}aussian, strongly dependent process, {\em Ann. Statist.\/}, 31, 2,
  586--612.

\bibitem{lmp:01}
Liseo B., Marinucci D. and Petrella L. (2001) Bayesian semiparametric inference
  on long-range dependence, {\em Biometrika\/}, 88, 4, 1089--1104.

\bibitem{lr:06}
Liseo B. and Rousseau J. (2006) Sequential importance sampling algorithm for
  {B}ayesian nonparametric long range inference, in: {\em Atti della XLIII
  Riunione Scientifica della Societ\`a Italiana di Statistica\/}, Societ\`a
  Italiana di Statistica, CLEUP, Padova, Italy, 43--46, vol. II.

\bibitem{mb:68}
Mandelbrot B.B. and Van~Ness J.W. (1968) Fractional {B}rownian motions,
  fractional noises and applications, {\em SIAM Rev.\/}, 10, 422--437.

\bibitem{mousou:03}
Moulines E. and Soulier P. (2003) Semiparametric spectral estimation for
  fractional processes, in: {\em Theory and applications of long-range
  dependence\/}, Birkh\"auser, Boston, MA, 251--301.

\bibitem{robi:91}
Robinson P.M. (1991) Nonparametric function estimation for long memory time
  series, in: {\em Nonparametric and semiparametric methods in econometrics and
  statistics (Durham, NC, 1988)\/}, Cambridge Univ. Press, Cambridge, Internat.
  Sympos. Econom. Theory Econometrics, 437--457.

\bibitem{robi:94}
Robinson P.M. (1994) Time series with strong dependence, in: {\em Advances in
  econometrics, Sixth World Congress, Vol.\ I (Barcelona, 1990)\/}, Cambridge
  Univ. Press, Cambridge, volume~23 of {\em Econom. Soc. Monogr.\/}, 47--95.

\bibitem{rob:95a}
Robinson P.M. (1995) Gaussian semiparametric estimation of long range
  dependence, {\em Ann. Statist.\/}, 23, 5, 1630--1661.

\bibitem{shw:01}
Shen X. and Wasserman L. (2001) Rates of convergence of posterior distibutions,
  {\em Annals of Statistics\/}, 29, 687--714.

\bibitem{whi:62}
Whittle P. (1962) Gaussian estimation in stationary time series, {\em Bull.
  Inst. Internat. Statist.\/}, 39, livraison 2, 105--129.

\end{thebibliography}

\section*{Acknowledgements}
Part of this work was done while the second Author was visiting the Universit\'e Paris Dauphine, CEREMADE. He thanks for warm hospitality and financial support.

\end{document}